\documentclass[11pt]{article}
\usepackage{amssymb,amsmath,amsthm}
\usepackage{epsfig}
\usepackage{verbatim}
\usepackage{graphicx}

\usepackage{color}
\usepackage{amssymb}

\newtheorem{thm}{Theorem}[section]

\newtheorem{prop}[thm]{Proposition}

\newtheorem{rem}[thm]{Remark}

\newcommand{\R}{\mathbb{R}}

\newcommand{\Pl}{\mathbb{P}}

\newcommand{\pv}{{\rm pv}}

\newcommand{\grad}{\nabla}

\newcommand{\dt}{\frac{d}{dt}}

\newcommand{\la}{\Lambda}

\newcommand{\te}{\theta}

\newcommand{\al}{\alpha}

\newcommand{\ep}{\varepsilon}
\newcommand{\gm}{\tilde{\gamma}}
\newcommand{\gn}{\gamma+\tilde{\varepsilon}}

\newcommand{\heatinv}{(\partial_t\!-\!\Delta)_0^{-1}}

\begin{document}

\title{
Regularity results for viscous 3D\\ Boussinesq temperature fronts}

\author{Francisco Gancedo and Eduardo Garc\'ia-Ju\'arez}

\date{\today}

\maketitle

\begin{abstract}
This paper is about the dynamics of non-diffusive temperature fronts evolving by the incompressible viscous Boussinesq system in $\R^3$. We provide local in time existence results for initial data of arbitrary size. Furthermore, we show global in time propagation of regularity for small initial data in critical spaces. The developed techniques allow to consider general fronts where the temperature is piecewise H\"older (not necessarily constant), which preserve their structure together with the regularity of the evolving interface.     
\end{abstract}

{\bf Keywords: }Boussinesq equations, temperature front, global regularity, singular heat kernels.

\setcounter{tocdepth}{1}

\section{Introduction}

In this paper we study the following active scalar equation
\begin{equation}\label{temperature}
\te_t + u\cdot\grad\te =0,
\end{equation}
where  $\theta(x,t)$ is the temperature of a three dimensional incompressible fluid
\begin{equation}\label{incompressible}
\grad\cdot u=0.
\end{equation}
The velocity $u(x,t)$ evolves by the viscous Boussinesq system
\begin{equation}\label{Boussinesq}
u_t +u\cdot\grad u-\Delta u+\nabla \pi=\te e_3,
\end{equation}
with $\pi(x,t)$ the fluid pressure and $e_3=(0,0,1)$. The viscosity and gravity constants are taken equal to one for the sake of simplicity and $(x,t)=(x_1,x_2,x_3,t)\in\R^3\times[0,+\infty)$.
The system above is the well-known Boussinesq equation \cite{Boussinesq1903} with viscosity and without heat diffusion. 

This system models natural convection phenomena generated by fluid flow due to the effect of buoyancy forces. Temperature gradients induce density variations from an equilibrium state, which gravity tends to restore. These flows are usually characterized by small deviations of the density with respect to a stratified reference state in hydrostatic balance. Potential energy is thus the main agent of movement, compared to inertia. Oberbeck was the first to notice by linearization that the buoyancy effect was proportional to temperature deviations \cite{Oberbeck1879}, and later Boussinesq \cite{Boussinesq1903} completed the model based on physical assumptions. It has been since then one of the main ingredients in geophysical models \cite{Gill82}, \cite{Majda03}, from ocean and atmosphere dynamics to mantle and solar inner convection, as well as a basic tool in building environmental engineering. In particular, it is an important model to understand the Rayleigh-B\'enard problem \cite{Constantin99}. 

From the mathematical point of view, the model is important due to the fact that it is related with the Euler and Navier-Stokes equations with constant density. Specifically the system (\ref{temperature},\ref{incompressible},\ref{Boussinesq}) contains 3D Navier-Stokes as a particular case and the inviscid 2D Boussinesq case corresponds to the 3D axisymmetric swirling Euler equations \cite{Majda02}. Furthermore, a very important feature is that, for the 2D and 3D cases in the Boussinesq system, vortex stretching mechanism is present. Therefore the well-posedness of the equations is a mayor open problem in the mathematical analysis of partial differential equations modeling incompressible fluids  \cite{Yudovich03}. 

For smooth initial data, the system was proved to exist for all time in the 2D case with regular initial data \cite{Chae06},\cite{Hou05}. The regularity of the initial data was improved later in \cite{Abidi07} and \cite{Hmidi07} using Besov and Sobolev spaces, respectively. Global-in-time results were shown through the scale of Sobolev spaces of different regularity in \cite{Hu15}. The uniqueness of weak solution was proved in \cite{Danchin08} making use of paradifferential calculus techniques.

On the other hand, in the 2D inviscid case the global existence of solution is still an open problem. There are numerical evidences of global-regularity \cite{E94} in the periodic setting but recent simulations indicated the possibility of finite-time blow-up in the case of bounded domain with regular boundary \cite{Luo14PNAS}. Based on this scenario, new one-dimensional Boussinesq models have been developed to show blow-up including boundary effects \cite{Choi15},\cite{Choi17}. They have been recently extended to dimension two \cite{Hoang16} for models which include the incompressibility condition \cite{Kiselev18}. Returning to the 2D inviscid Boussinesq system, it develops finite time blow-up in scenarios where the solutions have finite energy and evolve in spatial domains with a corner \cite{Elgindi17}. On the other hand, there is long-time existence for solutions  in the whole plane close to stable regimes, where the temperature of the fluid in increasing in the vertical direction \cite{Elgindi16}. Considering this setting for bounded domains the solutions have finite energy and the same result has been proven adding damping to the model \cite{Castro18}.

\par
There are also global existence results for the 2D Boussinesq system with anisotropic and partial viscosities (see \cite{Adhikari11},  \cite{Cao13}, \cite{Larios13}, \cite{Adhikari16}, \cite{Li16} and references therein). These scenarios model important physical situations with different horizontal and vertical scales for atmospheric and oceanic flows, where the viscosity constant can be zero in some axis-directions.

Taking the initial temperature trivial, $\theta(x,0)=0$, it is easy to obtain the Navier-Stokes equation from the Boussinesq approximation, so that in 3D the well-posedness is a challenging problem \cite{Fefferman06}. Similar blow-up criteria and global-in-time regularity for small initial data results can be shown \cite{Xu12}, \cite{Qiu10}, but the effect of gravity produces stratified temperature solutions in the inviscid case \cite{Widmayer15}.  

Here, we consider a free boundary problem governed by the incompressible Boussinesq system (\ref{temperature},\ref{incompressible},\ref{Boussinesq}), where the temperature has jumps of discontinuity.
 This setting provides important scenarios where the temperature is a front evolving with the fluid flow \cite{Gill82}, \cite{Majda03}. 
 In particular, our main concern is to study the propagation of regularity of the boundary of the front. There is a long tradition in the study of these kind of patch solutions, starting with the well-known vortex patch problem \cite{Chemin93}, \cite{Bertozzi93}. Moreover, having global existence for these scenarios in Boussinesq provides low-regular solutions for  Navier-Stokes  with an external force given by gravity, which are interesting by itself in the problems of global in time regularity. 
 
 This problem was first considered in \cite{Danchin17}, where initial fronts are given with regularity measured through the Besov spaces $\te_0\in B_{q,1}^{2/q-1}$, $q\in(1,2)$. The authors use paradifferential calculus and striated regularity techniques to obtain $C^{1+\gamma}$ propagation of regularity of the boundary of the fronts, $0<\gamma<1$, in 2D for arbitrary initial data and in 3D adding smallness assumptions. This result can then be applied to patch-type temperature, where the front takes different constant values on complementary domains. 

More generally, patch-type solutions have been highly studied for different equations coming from fluid mechanics models (see \cite{Chemin93}, \cite{Bertozzi93}, \cite{Cordoba10}, \cite{Fefferman16}, \cite{Gancedo14}, \cite{Gancedo18GarciaJuarez}).
In particular, in the setting of rapidly rotating temperature fronts, for the patch problem in 2D \cite{Rodrigo04} there is numerical evidence of pointwise collapse with curvature blow-up \cite{Cordoba05}. Moreover, it has been shown that the control of the curvature removes the possibility of pointwise interface collapse \cite{Gancedo14}.

In \cite{Gancedo17GarciaJuarez}, the authors proved that for 2D Boussinesq temperature patches the curvature of the interface cannot blow up in finite time. A new cancellation on the time-dependent singular integral operators given by the second derivatives of the solution to the heat equation was needed. A different proof for the persistence of low regularity was also shown and, by making use of level-set methods, an extra cancellation in the tangential direction is used to propagate higher regular interfaces.

In this paper we deal with 3D temperature fronts that do not need to be patches of constant temperature. We first provide a local-in-time existence result for very low regular initial temperature in Lebesgue spaces without size constraints (see theorem \ref{Case1} for more details). Adding a smallness condition in critical spaces, the results are global in time.
Notice that equation \eqref{Boussinesq} can be written as a forced heat equation and therefore the velocity is given by
\begin{equation}\label{decomposition}
u=e^{t\Delta}u_0-\heatinv\mathbb{P}(u\cdot\nabla u)+\heatinv\mathbb{P}(\theta e_3),
\end{equation}
where $\mathbb{P}$ is the Leray projector and $\heatinv f$ denotes the solution of the linear heat equation with force $f$ and zero initial condition:
\begin{equation*}
\heatinv f:=\int_0^t e^{(t-\tau)\Delta}f(\tau)d\tau.
\end{equation*}
Above we use the standard notation $e^{t\Delta}f=\mathcal{F}^{-1}(e^{-t|\xi|^2}\hat{f})$, where $\hat{}$ and $\mathcal{F}^{-1}$ denote Fourier transform and its inverse.
The particle trajectories of the system,
\begin{equation*}
\left\{\begin{aligned}
\frac{dX}{dt}(a,t)&=u(X(a,t),t),\\
X(a,0)&=a,
\end{aligned}\right.\mbox{ with the back-to-label map } A(X(a,t),t)=a,
\end{equation*}
gives the equation for the interface
$$Z_t(\alpha,t)=u(Z(\alpha,t),t),\quad \al\in\R^2.$$
Then, the regularity obtained for the velocity field (see Theorem \ref{Case1}) allows
to propagate the structure and interface regularity of fronts given by $\theta_0(x)=\theta_0(x)1_{D_0}(x)$, with $D_0\subset \R^3$ a bounded simply connected domain with boundary $\partial D_0\in C^{1+\gamma}$.
Under this well-posed scenario, we are able to prove that piecewise H\"older fronts,
\begin{equation*}
\theta_0(x)=\theta_1(x)1_{D_0}(x)+\theta_2(x)1_{D^c_0}(x),\hspace{0.5cm} \theta_1\in C^{\mu_1}(\overline{D}_0), \theta_2\in C^{\mu_2}(\overline{D^c_0})\cap L^1,\hspace{0.3cm}\mu_1,\mu_2\in (0,1), 
\end{equation*}
whose interface has bounded curvature, preserve this regularity locally in time for arbitrary data and globally in time for small data in critical spaces.
Taking advantage of the space-time singular integral given by the heat kernel, we find a new perspective to deal with the operators given by second derivatives of the temperature term in \eqref{decomposition}. This will allow to consider non-constant patches. The new viewpoint connects the \textit{parabolic} approach used in \cite{Gancedo17GarciaJuarez} with the \textit{elliptic} one in \cite{Gancedo18GarciaJuarez}, avoiding the use of time-weights and interpolation theory. Moreover, the technique provides a unified approach to propagate higher regular interfaces. Indeed,
for fronts with $C^{2+\gamma}$ boundary, where one cannot expect to gain enough regularity for the velocity globally in space, the approach used in 2D is no longer valid. That is due to the fact that in the 3D case the tangent vector fields are not divergence free. The \textit{parabolic-elliptic} approach overcomes this difficulty and shows in a clear way how to introduce contour dynamics methods to deal directly with the boundary evolution. Our new approach allows to use bootstrapping arguments getting propagation of regularity from weak solutions to $C^{1+\gamma}$, from $C^{1+\gamma}$ to $W^{2,\infty}$ and to higher regular interfaces.

\medskip
The paper is organized as follows: In Section \ref{intro}, we include the definition of the functional spaces used in the paper, a summary of paradifferential calculus and some regularity properties of the heat equation. In Section \ref{sec:2}, we give the local and global in time results for low regular temperature fronts, and show that $C^{1+\gamma}$ interfaces propagate preserving their regularity. Then, in Section \ref{sec:3}, we present the new approach that will allow us to deal with piecewise H\"older fronts whose initial curvature is bounded. We will find that $u\in L^1(0,T;W^{2,\infty})$ and therefore the boundary evolves without possibility of curvature blow up. In Section \ref{sec:4}, building upon the previous results, we will use contour dynamics methods to find the persistence of higher regularity.

\section{Functional Spaces and Preliminary Estimates}\label{intro}

We recall here the definition of Sobolev, H\"older and Besov spaces, together with some paradifferential calculus estimates and regularity properties of the heat equation  (see \cite{Bahouri11}, Chapters 1 and 2, for details).

Let $s\in\R$, $k$ a nonnegative integer and $\alpha\in(0,1)$. The homogeneous Sobolev space $\dot{H}^s(\R^3)$ is the space of tempered distributions $u$ with Fourier transform in $L^1_{loc}(\R^3)$ and such that
\begin{equation*}
\|u\|_{\dot{H}^s}^2=\int_{\R^3}|\xi|^{2s}|\hat{u}(\xi)|^2d\xi <\infty.
\end{equation*}
The nonhomogeneous counterpart is defined using the norm
\begin{equation*}
\|u\|_{H^s}=\|u\|_{L^2}+\|u\|_{\dot{H}^s}.
\end{equation*}
The H\"older space $C^{k+\alpha}$ is the space of $C^k$ functions $u$ such that 
\begin{equation*}
\|u\|_{C^{k+\alpha}}=\sup_{|\beta|\leq k}\|\partial^\beta u\|_{L^\infty}+\|u\|_{\dot{C}^{k+\alpha}}<\infty,
\end{equation*}
where the H\"older semi-norm $\|\cdot\|_{\dot{C}^{k+\alpha}}$ is defined by
\begin{equation*}
\|u\|_{\dot{C}^{k+\alpha}}=\sup_{|\beta|=k}\sup_{x\neq y}\frac{|\partial^\beta u(x)-\partial^\beta u(y)|}{|x-y|^\alpha}.
\end{equation*}
To define the homogeneous Besov spaces we first need to introduce the homogeneous Littlewood-Paley decomposition in $\R^3$.
Let $B=\{|\xi|\in\R^2: |\xi|\leq 4/3\}$ and $C=\{|\xi|\in\R^3:3/4\leq|\xi|\leq8/3\}$, and fix two smooth radial functions $\chi$ and $\varphi$ supported in $B$, $C$, respectively, and satisfying 
\begin{equation*}
\chi(\xi)+\sum_{j\geq0}\varphi(2^{-j}\xi)=1,\hspace{0.3cm}\forall \xi\in \R^3.
\end{equation*}
\begin{equation*}
\sum_{j\in\mathbb{Z}}\varphi(2^{-j}\xi)=1,\hspace{0.3cm}\forall \xi\in \R^3\setminus\{0\},
\end{equation*}
\begin{equation*}
|j-l|\geq2\rightarrow \text{Supp } \varphi(2^{-j}\cdot)\cap \text{Supp }\varphi(2^{-l}\cdot)=\emptyset.
\end{equation*}
The homogeneous dyadic blocks are defined by $\dot{\Delta}_jf=\mathcal{F}^{-1}(\varphi(2^{-j}\xi)\hat{f}(\xi))$.
Then, the homogeneous Besov space $\dot{B}_{p,q}^s(\R^3)$,  $p,q\in [1,\infty]$ consists of those tempered distributions $u\in S'_h(\R^3)$ for which
\begin{equation*}
\|u\|_{\dot{B}_{p,q}^s}=\|2^{js}\|\dot{\Delta}_j u\|_{L^p}\|_{l^q(j)}<\infty,
\end{equation*} 
where $u\in S'_h(\R^3)$ if $\lim_{\lambda\to 0}\|\mathcal{F}^{-1}(\theta(\lambda \xi)\hat{u}(\xi))\|_{L^\infty}=0$ for any $\theta$ smooth compactly supported function. 
\begin{rem}
It is well-known that the semi-norms $\|\cdot\|_{\dot{H}^s}$ and $\|\cdot\|_{\dot{B}_{2,2}^s}$ are equivalent, as well as  
$\|\cdot\|_{\dot{C}^{k+\alpha}}$ and $\|\cdot\|_{\dot{B}_{\infty,\infty}^{k+\alpha}}$.
\end{rem}
The next two propositions contains the embeddings that will be used along the paper.
\begin{prop}\label{lpembed}
	For any $p\in[1,\infty]$, $q_1\in [2,\infty)$, $q_2\in(1,2]$, the following continuous embeddings hold
\begin{equation*}
\begin{aligned}
\dot{B}_{p,1}^0\hookrightarrow L^p\hookrightarrow \dot{B}^0_{p,\infty},\qquad
\dot{B}^0_{q_1,2}\hookrightarrow L^{q_1},\qquad
L^{q_2}\hookrightarrow \dot{B}^0_{q_2,2}.
\end{aligned}
\end{equation*}	
Moreover, for $r_1\in[1,2]$, $r_2\in[2,\infty]$, it also holds that
\begin{equation*}
\dot{B}^0_{r_1,r_1}\hookrightarrow L^{r_1},\qquad  L^{r_2}\hookrightarrow \dot{B}^0_{r_2,r_2}.
\end{equation*}
\end{prop}
\begin{prop}\label{besovembed}
	Let $1\leq p_1\leq p_2\leq \infty$ and $1\leq q_1\leq q_2\leq \infty$. Then, for any $s\in\R$, the space $\dot{B}^s_{p_1,q_1}(\R^3)$ is continuously embedded in $\dot{B}^{s-3\big(\frac{1}{p_1}-\frac{1}{p_2}\big)}_{p_2,q_2}(\R^3)$.
\end{prop}
As particular cases, Proposition \ref{besovembed} together with Proposition \ref{lpembed} can be used to recover the standard Sobolev embeddings into H\"older and Lebesgue spaces.
 
We recall now some paradifferential calculus estimates. The homogeneous low-frequency cut-off operator $\dot{S}_j$ is defined by
\begin{equation*}
\dot{S}_ju=\mathcal{F}^{-1}(\chi(2^{-j}\xi)\hat{u}(\xi)).
\end{equation*} 
Considering then the Littlewood-Paley decompositions
\begin{equation*}
u=\sum_{l}\dot{\Delta}_lu,\quad v=\sum_{j}\dot{\Delta}_jv,
\end{equation*}
we can introduce Bony's decomposition for the product
\begin{equation*}
uv=\sum_{l,j}\dot{\Delta}_lu\dot{\Delta}_jv=\dot{T}_uv+\dot{T}_vu+\dot{R}(u,v),
\end{equation*}
where the homogeneous paraproduct $\dot{T}_uv$ and remainder $\dot{R}(u,v)$ of $u$ and $v$ are defined by
\begin{equation*}
\dot{T}_uv=\sum_j\dot{S}_{j-1}u\dot{\Delta}_jv,\quad \dot{R}(u,v)=\sum_{|l-j|\leq1}\dot{\Delta}_{l}u\dot{\Delta}_jv.
\end{equation*}
\begin{prop}\label{paradif}
For any $s, s_1,s_2\in \R$, $t<0$, $p,q\in[1,\infty]$, with
$$\frac{1}{q}=\frac{1}{q_1}+\frac{1}{q_2},\quad \frac{1}{p}=\frac{1}{p_1}+\frac{1}{p_2},$$
the following estimates hold
\begin{equation*}
\begin{aligned}
\|\dot{T}_uv\|_{\dot{B}^s_{p,q}}&\leq c \|u\|_{L^\infty}\|v\|_{\dot{B}^s_{p,q}},\\
\|\dot{T}_uv\|_{\dot{B}^{s+t}_{p,q}}&\leq c \|u\|_{\dot{B}^t_{\infty,q_1}}\|v\|_{\dot{B}^s_{p,q_2}},\\
\|\dot{R}(u,v)\|_{\dot{B}^{s_1+s_2}_{p,q}}&\leq c \|u\|_{\dot{B}^{s_1}_{p_1,q_1}}\|v\|_{\dot{B}^{s_2}_{p_2,q_2}}, \text{ if } s_1+s_2>0,\\
\|\dot{R}(u,v)\|_{\dot{B}^{s_1+s_2}_{p,\infty}}&\leq c \|u\|_{\dot{B}^{s_1}_{p_1,q_1}}\|v\|_{\dot{B}^{s_2}_{p_2,q_2}}, \text{ if } s_1+s_2\geq0 \text{ and } q=1.
\end{aligned}
\end{equation*}
\end{prop}
Using these estimates one can prove the following inequality in Sobolev spaces.
\begin{prop}
	For any $(s_1,s_2)\in (-3/2,3/2)$, $s_1+s_2\geq 0$, a constant $c$ exists such that
	\begin{equation}\label{SobolevParaDiff}
		\|uv\|_{\dot{H}^{s_1+s_2-3/2}}\leq c\|u\|_{\dot{H}^{s_1}}\|v\|_{\dot{H}^{s_2}}.
	\end{equation}
\end{prop}
Finally, we include some regularity estimates for the heat equation that will be used along the paper.
\begin{prop}\label{propoinitial} Let $s\geq0$, $r\in(1,\infty)$, $q\in\{1,\infty\}$, $\alpha\in(0,1)$ and $\varepsilon>0$. Then, the following estimates hold
	\begin{equation}\label{heatforcer}
		\|\nabla^2(\partial_t-\Delta)^{-1}_0f\|_{L^r_T(\dot{H}^{s})}\leq c\|f\|_{L^r_T(\dot{H}^s)},
	\end{equation}
	\begin{equation}\label{heatforces}	
		\|\nabla^2(\partial_t-\Delta)^{-1}_0f\|_{L^q_T(\dot{H}^s)}\leq c\|f\|_{L^q_T(\dot{H}^{s+\varepsilon})},
	\end{equation}
	\begin{equation}\label{heatforcesholder}	
\|\nabla^2(\partial_t-\Delta)^{-1}_0f\|_{L^1_T(\dot{C}^\alpha)}\leq c\|f\|_{L^1_T(\dot{C}^{\alpha+\varepsilon})},
\end{equation} 
	 	\begin{equation}\label{heatforcesbesov}	
	 \|\nabla^2(\partial_t-\Delta)^{-1}_0f\|_{L^\infty_T(\dot{B}^s_{2,\infty})}\leq c\|f\|_{L^\infty_T(\dot{B}^s_{2,\infty})},
	 \end{equation}
	\begin{equation}\label{heatinitial}
		\|\nabla^2 e^{t\Delta}u_0\|_{L^1_T(\dot{H}^s)}\leq c \|u_0\|_{\dot{H}^{s+\varepsilon}}.
	\end{equation}
	\begin{equation}\label{heatinitialholder}
	\|\nabla^2e^{t\Delta}u_0\|_{L^1_T(\dot{C}^\alpha)}\leq c \|u_0\|_{\dot{C}^{\alpha+\varepsilon}}.
	\end{equation}
	\begin{equation}\label{heatinitialbesov}
	\|e^{t\Delta}u_0\|_{L^\infty_T(\dot{B}^s_{2,\infty})}\leq c \|u_0\|_{\dot{B}^s_{2,\infty}}.
	\end{equation}
	Furthermore, there exists $u_0\in \dot{H}^s$ for which $\nabla^2e^{t\Delta}u_0\notin L^1_T(\dot{H}^s)$ and $u_0\in \dot{C}^\alpha$ for which $\nabla^2 e^{t\Delta}u_0\notin L^1_T(\dot{C}^\alpha).$

\end{prop}

Proof: The proof of \eqref{heatforcer} can be found in \cite{Krylov03}.  The proof of \eqref{heatforces}  follows from Bernstein inequalities and the decay of the heat kernel:
\begin{equation*}
	\begin{aligned}
		\|\nabla^2 \heatinv &f\|_{L^q_T(\dot{H}^s)}\leq \Big|\Big|\hspace{0.1cm} \Big|\Big| 2^{j(s+\varepsilon)}2^{j(2-\varepsilon)}c\int_0^t  e^{-c(t-\tau)2^{2j}}\|\Delta_j f\|_{L^2}(\tau)d\tau \Big|\Big|_{l^2(j)}\Big|\Big|_{L^q_T}\\
		&\leq  \Big|\Big| \int_0^t \frac{c}{(t-\tau)^{1-\varepsilon/2}}   \Big|\Big| 2^{j(s+\varepsilon)}\|\Delta_j f\|_{L^2}(\tau) \Big|\Big|_{l^2(j)}d\tau\Big|\Big|_{L^q_T}\leq c(T)\|f\|_{L^q_T(\dot{H}^{s+\varepsilon})}.
	\end{aligned}
\end{equation*}
In the last inequality above we have used Young's inequality for convolutions. The proof of \eqref{heatforcesholder} is similar and can be found in \cite{Gancedo17GarciaJuarez}.
On the other hand, there is no need of losing $\varepsilon$ derivatives if we work with Besov spaces with infinity as third index
\begin{equation*}
\begin{aligned}
\|\nabla^2 \heatinv f&\|_{L^\infty_T(\dot{B}^s_{2,\infty})}\leq \sup_{t\leq T}\sup_{j\in\mathbb{Z}} 2^{2j}c\int_0^t  e^{-c(t-\tau)2^{2j}}2^{sj}\|\Delta_j f\|_{L^2}(\tau)d\tau\\
&\leq c\sup_{j\in\mathbb{Z}} 2^{2j}\sup_{t\leq T} \int_0^t    e^{-c(t-\tau)2^{2j}}2^{sj}\|\Delta_j f\|_{L^2}(\tau) d\tau\leq c(T)\|f\|_{L^\infty_T(\dot{B}^{s}_{2,\infty})}.
\end{aligned}
\end{equation*}
We get \eqref{heatinitial} as before
\begin{equation*}
	\begin{aligned}
		\|\nabla^2 e^{t\Delta}u_0\|_{L^1_T(\dot{H}^s)}&\leq c \int_0^T \Big|\Big| 2^{j(s+\varepsilon)}2^{j(2-\varepsilon)}e^{-ct2^{2j}}\!\|\Delta_j u_0\|_{L^2}\Big|\Big|_{l_2(j)}dt\\
		&\leq  c\int_0^T \!\frac{\|u_0\|_{\dot{H}^{s+\varepsilon}}}{t^{1-\varepsilon/2}}dt\leq c(T)\|u_0\|_{\dot{H}^{s+\varepsilon}}.
	\end{aligned}
\end{equation*}
See \cite{Gancedo17GarciaJuarez} for \eqref{heatinitialholder}. The last estimate follows directly. The counterexamples in the last statements can be found in \cite{Fefferman17} and \cite{Gancedo17GarciaJuarez}, respectively.

\qed

\section{Local and global regularity results for $C^{1+\gamma}$ fronts}\label{sec:2}

This section is devoted to show a framework to provide local-in-time existence of low regular solutions for the Boussinesq system (\ref{temperature},\ref{incompressible},\ref{Boussinesq}) with no restriction on the size of the initial data. It also shows global-in-time solutions with smallness assumption on the initial data in critical spaces.

The proof is based on \text{a priori} energy estimates for weak solutions, then energy estimates for higher regularity and finally a bootstrapping argument using maximal regularity properties of the heat operator. For this last part we use the splitting \eqref{decomposition} as commented in the introduction.

\vspace{0.3cm}

\begin{thm}\label{Case1}
Assume $\gamma\in(0,1)$, $\ep\in(0,1-\gamma)$. Let $u_0\in H^{\frac12+\gamma+\ep}$ be a divergence-free vector field and $\te_0\in L^p$ for all $1\leq p<\frac{3}{1-\gamma-\varepsilon}$. Then, there is a unique solution $(u,\te)$ of (\ref{temperature},\ref{incompressible},\ref{Boussinesq}) with $u(x,0)=u_0$ such that
$$
\theta\in L^\infty(0,T;L^p),\quad u\in L^\infty(0,T;H^{\frac12+\gamma+\varepsilon})\cap L^2(0,T;H^{\frac32+\mu})\cap L^1(0,T;C^{1+\gamma+\tilde{\ep}}),
$$
for any $\mu\leq\min\{\gamma+\varepsilon,1/2\}$ and any $0<\tilde{\ep}<\ep$. 
The time of existence $T>0$ depends on the initial data in such a way that
$$
\int_{0}^T\| e^{\tau\Delta}u_0\|_{\dot{H}^1}^4(\tau)d\tau+\|\theta_0\|^2_{L^{3/2}}T<C_0,
$$
for $C_0>0$ an universal constant.

Furthermore, if the initial data satisfy
 \begin{equation}\label{smalldatahypo}
 \|u_0\|_{\dot{H}^{\frac12}}+\|\te_0\|_{L^1}<\delta,
 \end{equation}
for $\delta>0$ an universal constant, the solutions exist for all time $T>0$. 
\end{thm}

\begin{rem}\label{remark1}
The theorem above allows to propagate $C^{1+\gamma}$ regularity for temperature fronts; i.e. for initial $\theta_0(x)=\theta_0(x)1_{D_0}(x)$ with $D_0\subset \R^3$ a bounded simply connected domain with boundary $\partial D_0\in C^{1+\gamma}$,  $\theta_0\in L^p$, for all $1\leq p<\frac{3}{1-\gamma-\varepsilon}$, and $\partial D_0\in C^{1+\gamma}$. The temperature is given by
$$\te(x,t)=\te_0(A(x,t))1_{D(t)}(x) \hspace{0.2cm}{\rm{and}} \hspace{0.2cm} \partial D\in L^\infty(0,T;C^{1+\gamma}),$$
	where $D(t)=X(D_0,t)$.
\end{rem}

Proof: \\
\underline{\textbf{Local Existence}}: We consider first the $L^2$ energy balance for the Boussinesq system, obtaining that
\begin{equation*}
\frac12\dt\|u\|_{L^2}^2+\|\grad u\|_{L^2}^2=\int u_3\theta dx\leq \frac12\|\grad u\|^2_{L^2}+\frac12\|\te\|^2_{\dot{H}^{-1}}\leq \frac12\|\grad u\|^2_{L^2}+c\|\te\|^2_{L^{6/5}},
\end{equation*}
where the embedding $L^{6/5}\hookrightarrow \dot{H}^{-1}$ has been used.
From the transport character of \eqref{temperature} it is possible to find
\begin{equation*}
\|\te\|_{L^p}(t)\leq \|\te_0\|_{L^p},\quad \mbox{ for any }p\in [0,+\infty],
\end{equation*}
so that integration in time provides
\begin{equation}\label{l2balance}
\|u\|_{L^2}^2(t)+\int_0^T\|\grad u\|_{L^2}^2(\tau)d\tau\leq\|u_0\|^2_{L^2}+c\|\te_0\|^2_{L^{6/5}}T.
\end{equation}

Next we consider the $\dot{H}^{\frac12}$ norm evolution of the velocity. The argument is similar to \cite{Robinson16}, chapter 10, but the procedure is included for completeness.

Writing $u=v+w$ we decompose the velocity into a linear heat equation and a nonlinear system with zero initial data as follows
\begin{equation*}
v_t-\Delta v=0,\, v(x,0)=u_0(x);\qquad w_t-\Delta w=\Pl(u\cdot\nabla u)+\Pl(\theta e_3),\,w(x,0)=0,
\end{equation*}
where $\Pl$ denotes the Leray projection.
 It is then clear that
$$
\|v\|_{\dot{H}^{\frac12}}^2(t)+2\int_0^t\| v\|_{\dot{H}^{\frac32}}^2(\tau)d\tau=\|u_0\|_{\dot{H}^{\frac12}}^2.
$$
On the other hand
\begin{equation*}
\begin{aligned}
\frac12\frac{d}{dt}\|w\|_{\dot{H}^{\frac12}}^2+\|w\|_{\dot{H}^{\frac32}}^2&=\int \Lambda w\cdot (u\cdot\nabla u)dx +\int \la w_3\theta dx\\
&\leq \|w\|_{\dot{H}^{\frac32}}\|u\cdot\nabla u\|_{\dot{H}^{-\frac12}}+\| w\|_{\dot{H}^{\frac32}}\|\theta\|_{\dot{H}^{-\frac12}}.
\end{aligned}
\end{equation*}
The chain of bounds
$$
\|u\cdot\nabla u\|_{\dot{H}^{-\frac12}}\leq c\|u\cdot\nabla u\|_{L^\frac32}\leq c\|u\|_{L^6}\|\grad u\|_{L^2}\leq c\|u\|_{\dot{H}^{1}}^2,
$$
together with Young's inequality and the embedding $L^{3/2}\hookrightarrow \dot{H}^{-1/2}$ provide that
\begin{equation*}
\begin{aligned}
\frac{d}{dt}\|w\|_{\dot{H}^{\frac12}}^2+\|w\|_{\dot{H}^{\frac32}}^2&\leq  c\| u\|_{\dot{H}^{1}}^4+c_2\|\theta\|_{L^\frac32}^2\leq c_1\|w\|_{\dot{H}^{\frac12}}^2\|w\|_{\dot{H}^{\frac32}}^2+c_1\| v\|_{\dot{H}^{1}}^4+c_2\|\theta_0\|_{L^\frac32}^2.\\
\end{aligned}
\end{equation*}
The above inequality yields 
\begin{equation*}
\|w\|_{\dot{H}^{\frac12}}^2(t)+\int_0^t\|w\|_{\dot{H}^{\frac32}}^2(s)ds\leq \frac1{2c_1},\quad 0\leq t\leq T,
\end{equation*}
as long as
$$
\int_0^T(c_1\| v\|_{\dot{H}^{1}}^4(s)+c_2\|\theta_0\|_{L^\frac32}^2)ds<\frac1{4c_1}.
$$
In particular, for $t\in [0,T]$,
\begin{equation}\label{litb}
\|u\|_{\dot{H}^{\frac12}}^2(t)+\int_0^t\|u\|_{\dot{H}^{\frac32}}^2(s)ds+\int_0^t\|u\|_{\dot{H}^{1}}^4(s)ds\leq C(T).
\end{equation}
Next we consider the evolution of $1/2+\gm$ derivatives, with $\gm=\gamma+\ep$, as follows
\begin{equation*}
\begin{aligned}
\frac12\frac{d}{dt}\|u\|_{\dot{H}^{\frac12+\gm}}^2+\|u\|_{\dot{H}^{\frac32+\gm}}^2&\leq\int\la^{1+\gm}u\cdot\la^{\gm}(u\cdot \nabla u)dx+\int\la^{\frac32+\gm}u_3\la^{-\frac12+\gm}\theta dx\\
&\leq \|u\|_{\dot{H}^{1+\gm}}\|u\cdot\nabla u\|_{\dot{H}^{\gm}}+\|u\|_{\dot{H}^{\frac32+\gm}}\|\theta\|_{\dot{H}^{-\frac12+\gm}}.
\end{aligned}
\end{equation*}
Sobolev interpolation together with estimate \eqref{SobolevParaDiff} gives that
\begin{equation*}
\begin{aligned}
\frac12\frac{d}{dt}\|u\|_{\dot{H}^{\frac12+\gm}}^2+\|u\|_{\dot{H}^{\frac32+\gm}}^2&\leq \|u\|_{\dot{H}^{\frac12+\gm}}^\frac12\|u\|_{\dot{H}^{\frac32+\gm}}^\frac12\| u\|_{\dot{H}^{1}}\|u\|_{\dot{H}^{\frac32+\gm}}+\|u\|_{\dot{H}^{\frac32+\gm}}\|\theta\|_{\dot{H}^{-\frac12+\gm}} \\
&\leq \frac12\|u\|_{\dot{H}^{\frac32+\gm}}^2+\frac{3^3}4\|u\|_{\dot{H}^{\frac12+\gm}}^2\|u\|_{\dot{H}^{1}}^4+\|\theta\|_{\dot{H}^{-\frac12+\gm}}^2,
\end{aligned}
\end{equation*}
to find
\begin{equation}\label{1p2masgm}
\frac{d}{dt}\|u\|_{\dot{H}^{\frac12+\gm}}^2+\|u\|_{\dot{H}^{\frac32+\gm}}^2\leq \frac{3^3}2\|u\|_{\dot{H}^{\frac12+\gm}}^2\|u\|_{\dot{H}^{1}}^4+2\|\theta\|_{\dot{H}^{-\frac12+\gm}}^2.
\end{equation}
We consider to cases.\\

\underline{Case 1}: $\gm\in(0,\frac12]$. In this situation the bound $\|\theta\|_{\dot{H}^{-\frac12+\gm}}\leq c_3(p)\|\theta_0\|_{L^p}$ for $p=3/(2-\gm)$ together with \eqref{litb} allow us to use Gronwall's inequality in \eqref{1p2masgm} in order to conclude that
\begin{equation}\label{cl12masgm}
\|u\|_{\dot{H}^{\frac12+\gm}}^2(t)+\int_0^t\| u\|_{\dot{H}^{\frac32+\gm}}^2(s)ds\leq c_4(T),\qquad \forall t\in [0,T].
\end{equation}
Next we provide the solution of the system as follows
\begin{equation}\label{dnc}
 u= e^{t\Delta}u_0-\heatinv\Pl(u\cdot\nabla u)+\heatinv\Pl(\te e_3),
\end{equation}
to get for $\gamma<\gamma+\tilde{\varepsilon}<\gm=\gamma+\epsilon$ the following bound
\begin{equation*}
\begin{aligned}
\|u\|_{L^1_T(\dot{C}^{1+\gn})}&\leq \| e^{t\Delta}u_0\|_{L^1_T(\dot{C}^{1\!+\!\gn})}+\|\heatinv\Pl(u\cdot\nabla u)\|_{L^1_T(\dot{C}^{1+\gn})}\\
&\qquad +\|\heatinv\Pl(\te e_3)\|_{L^1_T(\dot{C}^{1+\gn})}.
\end{aligned}
\end{equation*}
Using Sobolev embedding and the fact that the Leray projector is bounded in H\"older spaces we can find
\begin{equation*}
\begin{aligned}
\|u\|_{L^1_T(\dot{C}^{1+\gn})} &\leq \|e^{t\Delta}u_0\|_{L^1_T(\dot{H}^{\frac52+\gn})}+\|\heatinv(u\cdot\nabla u)\|_{L^1_T(\dot{H}^{\frac52+\gn})}\\
&\qquad +\|\Delta\heatinv(\Delta^{-1}\te e_3)\|_{L^1_T(\dot{C}^{1+\gn})}.
\end{aligned}
\end{equation*}
Next we use \eqref{heatinitial}, \eqref{heatforces} and \eqref{heatforcesholder} to get
\begin{equation*}
\begin{aligned}
\|u\|_{L^1_T(\dot{C}^{1+\gn})} & \leq c(T)(\|u_0\|_{\dot{H}^{\frac12+\gm}}+\|u\cdot\nabla u\|_{L^1_T(\dot{H}^{\frac12+\gm})}+\|\te\|_{L^1_T(\dot{B}^{-1+\gm}_{\infty,\infty})}).
\end{aligned}
\end{equation*}
Using \eqref{SobolevParaDiff} and the embeddings $L^p\hookrightarrow\dot{B}^0_{p,p}\hookrightarrow\dot{B}^{-1+\tilde{\gamma}}_{\infty,\infty}$ from Propositions \ref{lpembed} and \ref{besovembed} with $p=\frac{3}{1-\gm}$, we obtain
\begin{equation*}
\begin{aligned}
\|u\|_{L^1_T(\dot{C}^{1+\gn})} & \leq c(T)(\|u_0\|_{\dot{H}^{\frac12+\gm}}+\|u\|^2_{L^2_T(\dot{H}^{\frac32+\gm})}+\|\te_0\|_{L^{\frac{3}{1-\gm}}})\leq c_4(T).
\end{aligned}
\end{equation*}
Therefore, we are done with the regularity for 
\begin{equation*}
u\in L^1_T(\dot{C}^{1+\gamma+\tilde{\varepsilon}}),\quad 0<\gamma<\gn<\gm\leq 1/2.
\end{equation*}
Using the same splitting above, an analogous computation  in nonhomogeneous spaces can be done to obtain that
\begin{equation}\label{l1tcmu}
u\in L^1_T(C^{1+\gamma+\tilde{\varepsilon}}),\quad 0<\gamma<\gn<\gm\leq 1/2.
\end{equation}

%

\underline{Case 2}: $\gm\in(\frac12,1)$. 
Using \eqref{l2balance} and \eqref{cl12masgm} with $\gm=1/2$ we find $u\in L^\infty_T(H^1)\cap L^2_T(H^2)$ so that interpolation provides
\begin{equation}\label{aux1}
u\in L^{\frac4{1+2\alpha}}_T(H^{\frac32+\alpha})\hookrightarrow L^{\frac4{1+2\alpha}}_T(C^{\alpha}),\quad 0<\alpha\leq 1/2,
\end{equation}
by Sobolev injection. Taking into account \eqref{aux1} with $\alpha=1/2$ and \eqref{l1tcmu} with $\gamma+\tilde{\varepsilon}=1/2-\bar{\varepsilon}\in (0,1/2)$, interpolation inequality
\begin{equation}
\|u\|_{C^\sigma}\leq c \|u\|_{C^\frac12}^{\lambda}\|u\|_{C^{\frac32-\bar{\ep}}}^{1-\lambda},\quad \frac12\leq\sigma<1, \, \quad\lambda=(\frac32-\sigma-\bar{\ep})/(1-\bar{\ep}), \quad 0<\bar{\varepsilon}< 1/2,
\end{equation}
provides that
\begin{equation*}
u\in L^p_T(C^\sigma), \hspace{0.5cm} p=\frac{4(1-\bar{\varepsilon})}{1+2(\sigma-\bar{\varepsilon})},\quad \frac12\leq\sigma<1, \quad 0<\bar{\varepsilon}< 1/2.
\end{equation*}
Therefore, by choosing $\bar{\varepsilon}=1-\tilde{\gamma}\in(0,1/2)$ and $\alpha=\frac{1-\tilde{\gamma}}{2\tilde{\gamma}}\in(0,1/2)$, we obtain that
\begin{equation*}
\begin{aligned}
\|u\otimes u\|_{L^1_T(\dot{C}^{\tilde{\gamma}})}&\leq 2\int_0^T \|u\|_{L^\infty}\|u\|_{\dot{C}^{\tilde{\gamma}}}dt\leq 2\int_0^T\|u\|_{C^\alpha}\|u\|_{C^{\tilde{\gamma}}}dt\\
&\leq 2\|u\|_{L^{4\tilde{\gamma}}_T\big(C^{\frac{1-\tilde{\gamma}}{2\tilde{\gamma}}}\big)}\|u\|_{L^{\frac{4\tilde{\gamma}}{4\tilde{\gamma}-1}}_T(C^{\tilde{\gamma}})}\leq c_4(T).
\end{aligned}
\end{equation*}
From \eqref{dnc} we find that
\begin{equation*}
\begin{aligned}
\|u\|_{L^1_T(\dot{C}^{1+\gn})}&\leq c(T)(\|u_0\|_{\dot{H}^{\frac12+\gm}}\!+\!\|\te_0\|_{L^\frac3{1-\gm}})\!+\!\|u\otimes u\|_{L^1_T(\dot{C}^{\tilde{\gamma}})}\leq c(T)(\|u_0\|_{\dot{H}^{\frac12+\gm}}\!+\!\|\te_0\|_{L^\frac3{1-\gm}})\!+\!c_4(T),
\end{aligned}
\end{equation*}
which, together with \eqref{l1tcmu}, yields that
\begin{equation*}
u\in L^1_T(C^{1+\gamma+\tilde{\varepsilon}}),\quad 0<\gamma<\gn<\gm=\gamma+\varepsilon<1.
\end{equation*}
Using \eqref{dnc} again, for $\varepsilon'\in (0,1-\gm)$,  we also obtain that
\begin{equation*}
\begin{aligned}
\|u\|_{L^\infty_T({\dot{H}}^{\frac12+\gm})}&\leq c(T)\left(
\|u_0\|_{\dot{H}^{\frac12+\gm}}+\|u\cdot\nabla u\|_{L^\infty_T(\dot{H}^{-\frac32+\gm+\varepsilon'})}+\|\te_0\|_{L^{\frac3{3-\gm-\varepsilon'}}}\right)\\
& \leq c(T)(\|u_0\|_{\dot{H}^{\frac12+\gm}}+\|u\|_{L^\infty_T(\dot{H}^{1})}\| u\|_{L^\infty_T(\dot{H}^{\gm+\varepsilon'})}+\|\te_0\|_{L^\frac3{3-\gm-\varepsilon'}})\leq c_4(T).
\end{aligned}
\end{equation*}
\\\underline{\textbf{Global Existence}}: We consider the splitting \eqref{dnc} and apply \eqref{heatinitialbesov} together with \eqref{heatforcesbesov} to find
\begin{equation*}
\begin{aligned}
\|u\|_{L_T^\infty(\dot{B}_{2,\infty}^{1/2})}&\leq c\|u_0\|_{\dot{B}_{2,\infty}^{1/2}}+k_1\|u\otimes u\|_{L_T^\infty(\dot{B}_{2,\infty}^{-1/2})}+c\|\te\|_{L_T^\infty(\dot{B}_{2,\infty}^{-3/2})}\\
&\leq c\|u_0\|_{\dot{B}_{2,\infty}^{1/2}}+k_1\|u\|^2_{L_T^\infty(\dot{B}_{2,\infty}^{1/2})}+c\|\te\|_{L_T^\infty(L^1)},
\end{aligned}
\end{equation*}
where we have used the paradifferential estimates of Proposition \ref{paradif} to bound the second term. The $L^p$ maximum principle for $\te$ together with the smallness condition
$$
\|u_0\|_{\dot{B}_{2,\infty}^{1/2}}+\|\te_0\|_{L^1}<\delta\leq \frac{1}{4k_1c}
$$ 
yield 
$$
\|u\|_{L_T^\infty(\dot{B}_{2,\infty}^{1/2})}\leq \frac{1}{2k_1}, \quad \forall\, T>0.
$$
The embedding $\dot{H}^\frac12\hookrightarrow \dot{B}_{2,\infty}^{1/2}$ allows to recover the more classical smallness assumption as state in the theorem.

Next we use the splitting \eqref{dnc} one more time, together with Proposition \ref{paradif} to obtain  
\begin{equation*}
\begin{aligned}
\|u\|_{L_T^2(\dot{H}^{\frac32})}&\leq c\|u_0\|_{\dot{H}^{\frac12}}+c\|u\otimes u\|_{L_T^2(\dot{H}^{\frac12})}+c\|\te\|_{L_T^2(\dot{H}^{-\frac12})}\\
&\leq c(\|u_0\|_{\dot{H}^{\frac12}}+T^{1/2}\|\te_0\|_{L^{3/2}})+c\|u\|_{L_T^\infty(\dot{B}_{2,\infty}^{1/2})}\|u\|_{L_T^2(\dot{H}^{\frac32})}.
\end{aligned}
\end{equation*}
Taking $k_1$ big enough it is possible to get
$$
\|u\|_{L_T^2(\dot{H}^{\frac32})}\leq C(\|u_0\|_{\dot{H}^{\frac12}}+T^{1/2}\|\te_0\|_{L^{3/2}}),\quad  \forall\, T>0.
$$
Finally, the bound $\|u\|_{\dot{H}^{1}}\leq c\|u\|^{1/2}_{\dot{B}^{1/2}_{2,\infty}}\|u\|^{1/2}_{\dot{H}^{\frac32}}$ gives
$$
\int_0^T\|u\|_{\dot{H}^{1}}^4(t)dt\leq \frac{C^2}{2k_1^2}(\|u_0\|_{\dot{H}^{\frac12}}^2+T\|\te_0\|_{L^{3/2}}^2),\quad  \forall\, T>0.
$$
It yields global existence in \eqref{1p2masgm} so that we can continue the proof in the same way as in the local-in-time approach.  

\hfill$\square$

\section{Local and global regularity results for $W^{2,\infty}$ fronts}\label{sec:3}
In this section we provide the local-in-time and global-in-time results to propagate the regularity of fronts with bounded curvature.

At this level of regularity, this problem can be considered as critical in the sense that one cannot expect more than $W^{2,\infty}$ regularity globally in space for the velocity since $\theta$ is merely bounded. Due to the singular integral operators given by two derivatives of the heat kernel, bounded functions would only yield $BMO$ type regular velocities, which are not generally bounded. Therefore, some extra cancellation is needed. 

The extra cancellation is achieved by the new elliptic-parabolic method that we introduced in this paper. 
Using the evolution equations together with integration by parts in time, we isolate the singularity of the space-time singular integrals. We reduce them to singular integrals only in space (fourth order Riesz transforms). These can be controlled thanks to the regularity provided by Theorem \ref{Case1} in the previous section, together with techniques for singular integrals with even kernels.

\begin{thm}\label{Case2}
	Let $u_0\in H^{\frac32+\ep}$ be a divergence-free vector field with $\ep\in(0,1)$. Assume that $D_0\subset \R^3$ is a bounded simply connected domain with boundary $\partial D_0\in W^{2,\infty}$, and $\te_0(x)=\te_0(x)1_{D_0}(x)$ with $\theta_0\in C^\mu(\overline{D}_0)$, $0<\mu<1$. 
	Then, there is a unique solution $(u,\te)$ of (\ref{temperature},\ref{incompressible},\ref{Boussinesq}) with $u(x,0)=u_0(x)$ such that
	$$\te(x,t)=\te_0(A(x,t))1_{D(t)}(x) \hspace{0.2cm}{\rm{and}} \hspace{0.2cm} \partial D\in L^\infty(0,T;W^{2,\infty}),$$
	where $D(t)=X(D_0,t)$. The regularity of the velocity is given by
	$$
	u\in L^\infty(0,T;H^{\frac32+\varepsilon})\cap L^1(0,T;W^{2,\infty}).$$ 
	The time of existence $T>0$ depends on the initial data in such a way that
	$$
	\int_{0}^T\| e^{\tau\Delta}u_0\|_{\dot{H}^1}^4(\tau)d\tau+\|\theta_0\|^2_{L^{3/2}}T<C_0,
	$$
	for $C_0>0$ an universal constant. Furthermore, if the initial data satisfy
	\begin{equation*}
	\|u_0\|_{\dot{H}^{\frac12}}+\|\te_0\|_{L^1}<\delta,
	\end{equation*}
	for $\delta>0$ an universal constant, the solutions exist for all time $T>0$. 
\end{thm}

\begin{rem}\label{remark42}
	The theorem above allows to show an analogous result for initial fronts of the form
	$\theta_0(x)=\theta_1(x)1_{D_0}(x)+\theta_2(x)1_{D^c_0}(x)$ with $\theta_1\in C^{\mu_1}(\overline{D}_0)$, $\theta_2\in C^{\mu_2}(\overline{D^c_0})\cap L^1$ and $\mu_1,\mu_2\in (0,1)$. Then, the same conclusions for $u$ and $\theta$ are obtained and the front propagates as
	\begin{equation*}
	\te(x,t)=\te_1(A(x,t))1_{D(t)}(x)+\theta_2(A(x,t))1_{D^c(t)}(x) \hspace{0.2cm}{\rm{with}} \hspace{0.2cm} \partial D\in L^\infty(0,T;W^{2,\infty}).
	\end{equation*}
\end{rem}

Proof: First point in the argument is to use Theorem \ref{Case1} in order to find a unique solution up to a time $T>0$ to the system. Under the smallness assumption the previous estimates are global and so are the following, giving the global existence result. 

Next we use splitting \eqref{dnc} to find
\begin{equation}\label{aux}
\|u\|_{L^\infty_T(\dot{H}^{\frac32+\varepsilon})}\leq \|u_0\|_{\dot{H}^{\frac32+\varepsilon}}+c\|u\cdot\nabla u\|_{L^\infty_T(\dot{H}^{-\frac12+\varepsilon'})}+c\|\te\|_{L^\infty_T(\dot{H}^{-\frac12+\varepsilon'})},
\end{equation}
where $0<\varepsilon<\varepsilon'<1$.
For $0<\varepsilon'\leq 1/2$, the last term is bounded by Sobolev embedding $L^{3/(2-\varepsilon')}\hookrightarrow \dot{H}^{-\frac12+\varepsilon'}$. 
For $1/2<\varepsilon'<1$, we use that Theorem \ref{Case1} guarantees that $\theta(t)$ is a H\"older patch with $C^{1+\gamma}$ boundary for any $0<\gamma<1$. More specifically, we recall that $\theta(x,t)=\theta_0(A(x,t))1_{D(t)}(x)$, with $\theta_0(A(x,t))\in C^\mu(\overline{D}(t))$. In particular, $\theta_0(A(x,t))\in H^\mu(D(t))$ and therefore there is an extension of it to $\mathbb{R}^3$  (see \cite{McLean00}), which we denote  $\tilde{\theta}(t)\in H^\mu\cap L^\infty$.
Then, for $0<s<\min\{\mu,\frac12\}$, we make use of standard paradifferential calculus estimates (see Proposition \ref{paradif}) to find that
\begin{equation*}
	\begin{aligned}
		\|\theta(t)\|_{\dot{H}^{s}}&= \|T_{1_{D(t)}}\tilde{\theta}(t)\|_{\dot{H}^s}+\|T_{\tilde{\theta}(t)} 1_{D(t)}\|_{\dot{H}^s}+\|R(1_{D(t)},\tilde{\theta}(t))\|_{\dot{H}^s}\\
		&\leq c(\|1_{D(t)}\|_{L^\infty}\|\tilde{\theta}\|_{\dot{H}^s}+\|\tilde{\theta}(t)\|_{L^\infty}\|1_{D(t)}\|_{\dot{H}^s})\leq C(T),
	\end{aligned}
\end{equation*}
where in the last step we have use the fact that the characteristic function of a bounded Lipschitz domain is in $\dot{H}^s$ for $0< s<\frac12$ (see \cite{Faraco13}). Then,  using \eqref{SobolevParaDiff} for the second term in \eqref{aux} we obtain
\begin{equation*}
\begin{aligned}
\|u\|_{L^\infty(\dot{H}^{\frac32+\varepsilon})}&\leq \|u_0\|_{\dot{H}^{\frac32+\varepsilon}}+c\|u\otimes u\|_{L^\infty_T(\dot{H}^{\frac12+\varepsilon})}+C(T)\\
&\leq \|u_0\|_{\dot{H}^{\frac32+\varepsilon}}+c\|u\|^2_{L^\infty_T(\dot{H}^{1+\frac{\varepsilon}{2}})}+C(T)\leq C(T),
\end{aligned}
\end{equation*}
where the last terms are controlled due to the estimates found in the previous section.

The last estimate for the velocity is performed:
\begin{equation}\label{W2inf}
\begin{aligned}
\|u\|_{L^1_T(\dot{W}^{2,\infty})}&\leq c( \|u_0\|_{H^{\frac32+\ep}}+\|u\cdot\nabla u\|_{L^1_T(C^{\ep})}+\|\nabla^2\heatinv\Pl(\te e_3)\|_{L^1_T(L^\infty)})\\
&\leq c(\|u_0\|_{H^{\frac32+\ep}}+\|u\|_{L^\infty_T(C^{\ep})}\|\nabla u\|_{L^1_T(C^{\ep})}+\|\nabla^2\heatinv\Pl(\te e_3)\|_{L^1_T(L^\infty)})\\
&\leq C(T)+\|\nabla^2\heatinv\Pl(\te e_3)\|_{L^1_T(L^\infty)},
\end{aligned}
\end{equation}
so that it remains to bound the last term above.

Next we analyze the singular integral operator $\nabla^2\heatinv\Pl(\theta e_3)$. We apply the Fourier transform to find 
$$
\mathcal{F}(\nabla^2\heatinv\Pl(\theta e_3))(\xi,t)=\xi_j\xi_k\Big(\delta_{l,3}-\frac{\xi_l\xi_3}{|\xi|^2}\Big)\int_0^t\exp(-|\xi|^2(t-\tau))\hat{\theta}(\xi,\tau)d\tau,
$$
for $j$, $k$, $l$ form 1 to 3 and $\delta_{l,3}$ the Kronecker delta. It shows that we only need to deal with the following four cases 
$$
\partial_1^3\partial_3(-\Delta)^{-1}\heatinv\theta(x,t)=\pv\int_0^t\int_{\R^3}K_1(x-y,t-\tau)\theta(y,\tau)dyd\tau,
$$
$$ 
\partial_1^2\partial_2\partial_3(-\Delta)^{-1}\heatinv\theta(x,t)=\pv\int_0^t\int_{\R^3}K_2(x-y,t-\tau)\theta(y,\tau)dyd\tau,
$$
$$
\partial_1^2\partial_3^2(-\Delta)^{-1}\heatinv\theta(x,t)=\pv\int_0^t\int_{\R^3}K_3(x-y,t-\tau)\theta(y,\tau)dyd\tau,
$$

$$
\partial_1^4(-\Delta)^{-1}\heatinv\theta(x,t)=\pv\int_0^t\int_{\R^3}K_4(x-y,t-\tau)\theta(y,\tau)dyd\tau,
$$
by exchanging coordinates.
The principal values above are understood as a limit removing the time singularity. The identity $$\Delta K=\partial_t K,$$ for $K$ the heat kernel at any $t\neq 0$, 
implies that
\begin{equation}\label{kernels}
\begin{aligned}
\widehat{K_1}(\xi,t)&=\frac{\xi_1^3\xi_3}{|\xi|^4}\widehat{\partial_t K}(\xi,t), \qquad \widehat{K_2}(\xi,t)= \frac{\xi_1^2\xi_2\xi_3}{|\xi|^4}\widehat{\partial_tK}(\xi,t),\\
\widehat{K_3}(\xi,t)&=\frac{\xi_1^2\xi_3^2}{|\xi|^4}\widehat{\partial_t K}(\xi,t), \qquad \widehat{K_4}(\xi,t)=\frac{\xi_1^4}{|\xi|^4}\widehat{\partial_t K}(\xi,t).
\end{aligned}
\end{equation}
Therefore, we can write (see Chapter 3.3 in \cite{Stein70})
\begin{equation}\label{Ki}
K_i(x,t)=\partial_t(k_i(x)* K(x,t))+\frac1{15}\delta_{i3}\partial_t K(x,t)+\frac1{5}\delta_{i4}\partial_t K(x,t)
\end{equation}
where the kernels $k_i$, $i=1,..., 4$, correspond to fourth order Riesz transforms, hence they are even, homogeneous of degree $-3$, and have zero mean on spheres. We denote by $\delta_{ij}$ the Kronecker delta.

Going back to \eqref{W2inf}, we can now bound the temperature terms by the following
\begin{equation}\label{i1234}
\|\nabla^2\heatinv\Pl(\te e_3)\|_{L^1_T(L^\infty)}\leq c \sum_{i=1}^4\left|\left|\pv \int_0^t \int_{\mathbb{R}^3} K_i(x-y,t-\tau)\theta(y,\tau)dyd\tau\right|\right|_{L^1_T(L^\infty)}.
\end{equation}
We show the details in the case $i=4$, as the rest can be handled in a similar manner. 

Then we start by splitting as follows
\begin{equation*}
\pv \int_0^t \int_{\mathbb{R}^3} K_4(x-y,t-\tau)\theta(y,\tau)dyd\tau=I_1+I_2,
\end{equation*}
where
\begin{equation}\label{I2}
\begin{aligned}
I_1&=-\pv \int_0^t \int_{\mathbb{R}^3} \partial_\tau ((k_4* K)(x-y,t-\tau))\theta(y,\tau)dyd\tau,\\
I_2&=-\frac15\pv \int_0^t \int_{\mathbb{R}^3} \partial_\tau K(x-y,t-\tau))\theta(y,\tau)dyd\tau.
\end{aligned}
\end{equation}
Using the equation \eqref{temperature}, the term $I_1$ becomes
\begin{equation*}
\begin{aligned}
I_1&=-\pv \int_0^t \int_{\mathbb{R}^3} \partial_\tau ((k_4* K)(x-y,t-\tau)\theta(y,\tau))dyd\tau\\
&\quad-\pv\int_0^t \int_{\mathbb{R}^3} (k_4* K)(x-y,t-\tau)\nabla\cdot(u(y,\tau)\theta(y,\tau)) dyd\tau,
\end{aligned}
\end{equation*}
so integration by parts shows that
\begin{equation}\label{I1}
\begin{aligned}
I_1&=-\lim_{\tau\to t^-}\int_{\mathbb{R}^3}(k_4* K)(x-y,t-\tau)\theta(y,\tau)dy+\int_{\mathbb{R}^3}(k_4* K)(x-y,t)\theta(y,0)dy\\
&\quad+\int_0^t\int_{\mathbb{R}^3}\nabla K(x-y,t-\tau)\cdot (k_4* (u \theta))(y,\tau)dyd\tau\\
&=J_1+J_2+J_3.
\end{aligned}
\end{equation}
The first term can be written as follows
\begin{equation*}
\begin{aligned}
J_1&=-\lim_{\tau\to t^-}\int_{\mathbb{R}^3}K(x-y,t-\tau) (k_4* \theta)(y,\tau)dy=-\lim_{\epsilon \to 0^+} K(\epsilon)* (k_4* \theta(t-\epsilon))(x)\\
&=-(
k_4* \theta) (x,t).
\end{aligned}
\end{equation*}
We note here that $k_4$ defines a singular integral operator, and thus $J_1$ is not bounded for a general bounded function. We first define a cut-off distance
\begin{equation}\label{cutoff}
\delta=\min_{\tau\in[0,t]}\left(\frac{|\nabla Z|_{\inf}(\tau)}{\|\nabla Z\|_{C^\eta}(\tau)}\right)^{1/\eta},\hspace{1cm}\eta\in(0,1),
\end{equation}
where 
$$
|\nabla Z|_{\inf}(t)=\min_{j}\Big\{\min\Big\{\inf_{\al\in \mathcal{N}_j} |\partial_{\al_1}Z(\al,t)|, \inf_{\al\in \mathcal{N}_j} |\partial_{\al_2}Z(\al,t)|,\inf_{\al,\beta\in \mathcal{N}_j,\al\neq\beta} \frac{|Z(\al,t)-Z(\beta,t)|}{|\al-\beta|}\Big\}\Big\}.
$$
Above the neighborhoods $\mathcal{N}_j$, $j=1,...,L$, provide local charts of the free boundary $\partial D(t)$ so that for any $x\in\partial D(t)$ there exists a $\mathcal{N}_j\subset\R^2$ such that $x=Z(\al,t)$ with $\al\in\mathcal{N}_j$. The positive quantity $\delta$ is fixed due to Theorem \ref{Case1}. Then we can perform the following splitting
\begin{equation}\label{J1}
J_1=-\int_{D(t)\cap\{|x-y|\geq \delta \}}k_4(x-y)\theta(y,t)dy-\pv\int_{ D(t)\cap \{|x-y|< \delta \}}k_4(x-y)\theta(y,t)dy=L_1+L_2,
\end{equation}
where the first term is bounded by
\begin{equation}\label{L1}
|L_1|\leq \int_{|x-y|\geq \delta}|k_4(x-y)||\theta(y,t)|dy\leq c \|\theta_0\|_{L^\infty}|\log{\delta}||D_0|.
\end{equation}
In order to bound $L_2$ we distinguish between two cases: $x\in \overline{D(t)}$ and $x\notin \overline{D(t)}$. From now on, if $x\in\partial D(t)$ the meaning of $\theta(x,t)$ is the limit of $\theta(y,t)$ from inside $D(t)$ as $y\to x$.

In the first case, we split $L_2$ as follows
\begin{equation*}
L_2=-\pv\int_{D(t)\cap\{|x-y|< \delta \}}k_4(x-y) (\theta(y,t)-\theta(x,t))dy-\theta(x,t)\pv\int_{D(t)\cap\{|x-y|< \delta \}} k_4(x-y)dy,
\end{equation*}
and therefore
\begin{equation*}
|L_2|\leq c \|\theta(t)\|_{C^\mu(\overline{D(t)})}\frac{\delta^\mu}{\mu}+\Big|\theta(x,t)\pv\int_{D(t)\cap\{|x-y|< \delta \}} k_4(x-y)dy\Big|.
\end{equation*}
Since $\theta$ satisfies a transport equation, the regularity obtained for $u$ allows us to find that
\begin{equation*}
\|\theta(t)\|_{C^\mu(\overline{D(t)})}\leq \|\theta_0\|_{C^\mu(\overline{D}_0)}e^{c\int_0^t\|\nabla u\|_{L^\infty}(\tau)d\tau}\leq c(\|u_0\|_{H^{\frac12+\varepsilon}},T)\|\theta_0\|_{C^\mu(\overline{D}_0)},
\end{equation*}
so we have that
\begin{equation}\label{L2}
|L_2|\leq c(\|u_0\|_{H^{\frac12+\varepsilon}},T)\|\theta_0\|_{C^\mu(\overline{D}_0)}+\Big|\theta(x,t)\pv\int_{D(t)\cap\{|x-y|< \delta \}} k_4(x-y)dy\Big|.
\end{equation}
The last term above can be bounded since the kernels are homogeneous and even, and $\partial D(t)\in C^{1+\gamma}$ due to Remark \ref{remark1} (see \cite{Bertozzi93} for the complete argument in 2d, and \cite{Cordoba10} for its extension to the three dimensional case).

In the case $x\notin \overline{D(t)}$, we define $\tilde{x}=\arg d(x, \partial D(t))\in \partial D(t)$. Then we can split $L_2$ as follows
\begin{equation*}
L_2=-\pv\int_{D(t)\cap\{|x-y|< \delta \}}k_4(x-y) (\theta(y,t)-\theta(\tilde{x},t))dy-\theta(\tilde{x},t)\pv\int_{D(t)\cap\{|x-y|< \delta \}} k_4(x-y)dy,
\end{equation*}
so taking into account the triangle inequality we find that
\begin{equation*}
\begin{aligned}
|L_2|&\leq \|\theta(t)\|_{C^\mu(\overline{D(t)})}\int_{D(t)\cap\{|x-y|< \delta \}}|k_4(x-y)|2^\mu|x-y|^{\mu}dy\\
&\quad+\Big|\theta(\tilde{x},t)\pv\int_{D(t)\cap \{|x-y|< \delta \}}k_4(x-y)dy\Big|,
\end{aligned}
\end{equation*}
thus we conclude that
\begin{equation}\label{L2bound}
|L_2|\leq c(\|u_0\|_{H^{\frac12+\varepsilon}},\|\theta_0\|_{C^\mu(\overline{D}_0)},\delta,T).
\end{equation}
Going back to \eqref{J1}, the bounds \eqref{L1} and \eqref{L2bound} give
\begin{equation}\label{J1bound}
\|J_1\|_{L^\infty_T(L^\infty)}\leq c(\|u_0\|_{H^{\frac12+\varepsilon}},\|\theta_0\|_{C^s(\overline{D}_0)},\delta,|D_0|, T).
\end{equation}
We proceed now to bound $J_2$ in \eqref{I1}. We can write it in the following manner
\begin{equation*}
J_2=\int_{\mathbb{R}^3} K(x-y,t) (k_4*\theta_0)(y)dy=e^{t\Delta}(k_4*\theta_0)(x),
\end{equation*}
so standard properties of the heat equation give us that
\begin{equation*}
\|J_2\|_{L^\infty_T(L^\infty)} \leq\|k_4*\theta_0\|_{L^\infty}.
\end{equation*}
Using the same reasoning as in the term $J_1$, we can conclude that
\begin{equation}\label{J2bound}
\|J_2\|_{L^\infty_T({L^\infty})}\leq c(\|\theta_0\|_{C^s(\overline{D}_0)},\delta,|D_0|).
\end{equation}
The last term is indeed more regular and it can be bounded as follows
\begin{equation}
\begin{aligned}\label{J3bound}
\|J_3\|_{L^\infty_T(L^\infty)}&\leq \Big|\Big|\int_0^t \|\nabla K(t-\tau)\|_{L^{\frac{3}{2+\varepsilon}}}\|k_4*(u\theta)(\tau)\|_{L^{\frac{3}{1-\varepsilon}}}d\tau\Big|\Big|_{L^\infty_T}\leq c(T)\|u\theta\|_{L^\infty_T(L^{\frac{3}{1-\varepsilon}})}\\
&\leq c(\|\theta_0\|_{L^\infty},\|u_0\|_{H^{\frac12+\varepsilon}},T).
\end{aligned}
\end{equation}
So introducing the bounds \eqref{J1bound}, \eqref{J2bound} and \eqref{J3bound} in \eqref{I1} we have that
\begin{equation*}
\|I_1\|_{L^\infty_T(L^\infty)}\leq c(\|u_0\|_{H^{\frac12+\varepsilon}},\|\theta_0\|_{C^s(\overline{D}_0)},\delta,|D_0|, T).
\end{equation*}
The term $I_2$ is analogous to $I_1$ but replacing the Riesz transforms with identities, so we find that
\begin{equation*}
\|I_2\|_{L^\infty_T(L^\infty)}\leq c(\|\theta_0\|_{L^\infty},\|u_0\|_{H^{\frac12}},T).
\end{equation*}
We are then done with estimate \eqref{i1234} and therefore plugging it into \eqref{W2inf} we find the regularity needed to end the proof. 

\qed

\section{Local and global regularity results for $C^{2+\gamma}$ fronts}\label{sec:4}

This section is devoted to prove the results for $C^{2+\gamma}$ fronts.
As commented in the beginning of the previous section, one cannot expect to obtain $C^{2+\gamma}$ regularity globally in space for the velocity, since $\theta$ is only bounded.
Indeed, taking two derivatives in \eqref{dnc}, the hardest part is to study the regularity of the last term on  the boundary of the H\"older-patch.
We first use the trick of the previous section to reduce the evolution terms of the velocity into singular integrals of Riesz type at the fixed time $t$. Then, we study the H\"older regularity of these singular integrals on the boundary. After some technical splittings, we take advantage of the fact that fourth order Riesz transform kernels can be integrated (see $L_6$ term, \eqref{L6}) to fully introduce a contour dynamics formulation.
Studying these new kernels together with the previous regularity results, the new theorem below follows.

\begin{thm}\label{Case3}
	Let $u_0\in H^{\frac32+\gamma+\ep}$ be a divergence-free vector field with $\gamma\in(0,1)$ and $0<\ep<\min\{1/2,1-\gamma\}$. Assume that $D_0\subset \R^3$ is a bounded simply connected domain with boundary $\partial D_0\in C^{2+\gamma}$, and $\te_0(x)=\te_0(x)1_{D_0}(x)$ with $\theta_0\in C^\gamma(\overline{D}_0)$. 
	Then, there is a unique solution $(u,\te)$ of (\ref{temperature},\ref{incompressible},\ref{Boussinesq}) with $u(x,0)=u_0(x)$ such that
	$$\te(x,t)=\te_0(A(x,t))1_{D(t)}(x) \hspace{0.2cm}{\rm{and}} \hspace{0.2cm} \partial D\in L^\infty(0,T;C^{2+\gamma}),$$
	where $D(t)=X(D_0,t)$. The regularity of the velocity is given by
	$$
	u\in L^\infty(0,T;H^{\frac32+\gamma+\varepsilon})\cap L^1(0,T;W^{2,\infty}).$$
	The time of existence $T>0$ depends on the initial data in such a way that
	$$
	\int_{0}^T\| e^{\tau\Delta}u_0\|_{\dot{H}^1}^4(\tau)d\tau+\|\theta_0\|^2_{L^{3/2}}T<C_0,
	$$
	for $C_0>0$ an universal constant. Furthermore, if the initial data satisfy
	\begin{equation*}
	\|u_0\|_{H^{\frac12}}+\|\te_0\|_{L^1}<\delta,
	\end{equation*}
	for $\delta>0$ an universal constant, the solutions exist for all time $T>0$. 
\end{thm}

\begin{rem}
	The theorem above allows to show an analogous result for initial fronts of the form
	$\theta_0(x)=\theta_1(x)1_{D_0}(x)+\theta_2(x)1_{D^c_0}(x)$ with $\theta_1\in C^{\gamma}(\overline{D}_0)$, $\theta_2\in C^{\gamma}(\overline{D^c_0})\cap L^1$. The same  conclusions for $u$ and $\theta$ are obtained and the front propagates as 
	$$\te(x,t)=\te_1(A(x,t))1_{D(t)}(x)+\theta_2(A(x,t))1_{D^c(t)}(x) \hspace{0.2cm}{\rm{with}} \hspace{0.2cm} \partial D\in L^\infty(0,T;C^{2+\gamma}).$$
\end{rem}

Proof: The regularity $u\in L^\infty_T(H^{\frac32+\gamma+\varepsilon})$ follows from Theorem \ref{Case2}. We need to study the $C^{2+\gamma}$ regularity of the velocity. From the equation, since $\theta$ is not continuous, one cannot expect to obtain such regularity globally in space, so we need to study it on the surface. The nonlinear term and the corresponding to the initial data in the splitting \eqref{dnc} can be treated as in \eqref{aux}. In particular, for these terms one can indeed obtain the regularity globally in space:
\begin{equation*}
\begin{aligned}
\|e^{t\Delta}u_0\|_{L^1_T(C^{2+\gamma})}&\leq \|u_0\|_{H^{\frac32+\gamma+\varepsilon}},\\
\|\heatinv\mathbb{P}(u\cdot\nabla u)\|_{L^1_T(C^{2+\gamma})}&\leq c\|u\cdot\nabla u\|_{L^1_T(C^{\gamma+\varepsilon})}\leq c\|u\|_{L^\infty_T(C^{\gamma+\varepsilon})}\|\nabla u\|_{L^1_T(C^{\gamma+\varepsilon})}\leq C(T).
\end{aligned}
\end{equation*}
It remains to deal with the temperature term. In order to deal with it, we consider as before one of the main kernels as explained in the proof of Theorem \ref{Case2}. The others can be treated in a similar manner. We deal with the kernel $K_4$ given in \eqref{kernels}. Integration by parts in time provides as before 
\begin{equation}\label{I1I2I3decomp}
\pv \int_0^t \int_{\mathbb{R}^3} K_4(x-y,t-\tau)\theta(y,\tau)dyd\tau=I_1+I_2+I_3,
\end{equation}
where
\begin{equation}\label{I1I2I3}
\begin{aligned}
I_1&=-(k_4* \theta) (x,t)-\frac1{5}\theta(x,t),\\
I_2&=e^{t\Delta}(k_4*\theta_0)(x)+\frac1{5}e^{t\Delta}\theta_0(x),\\
I_3&=\int_0^t\int_{\mathbb{R}^3}\nabla K(x-y,t-\tau)\cdot ((k_4+\frac{1}{5}\delta_0)* (u \theta))(y,\tau)dyd\tau,
\end{aligned}
\end{equation}
with $k_4$ given in \eqref{Ki} and $\delta_0$ the Dirac delta.
Notice that the second and third term correspond to solutions of the linear heat equation, and therefore can be bounded as follows
\begin{equation}\label{I2I3bound}
\begin{aligned}
\|I_2\|_{L^1_T(C^\gamma)}&\leq c(T)\|\theta_0\|_{L^\infty},\\
\|I_3\|_{L^1_T(C^\gamma)}&\leq c(T)\|u\|_{L^1_T(L^\infty)}\|\theta\|_{L^\infty_T(L^\infty)}\leq c(T,\|u_0\|_{H^{\frac12+\varepsilon}},\|\theta_0\|_{L^\infty}).
\end{aligned} 
\end{equation}


Therefore, it only remains to deal with the term $I_1$. Since we want to study the H\"older regularity along the surface, we consider two points on the surface $x=Z(\alpha,t)$, $x+h=Z(\tilde{\alpha},t)$. Then we start with the following splitting to deal with the H\"older norm
\begin{equation*}
I_1(x+h)-I_1(x)=J_1+J_2,
\end{equation*}
where	
\begin{equation*}
J_1=\int_{D(t)}(k_4(x-y)-k_4(x+h-y))\theta(y,t)dy,
\end{equation*}
and
\begin{equation*}
J_2=\frac{1}{5}(\theta(x)-\theta(x+h)).
\end{equation*}
The term above gives
$$
|J_2|\leq \frac{1}{5}\|\theta\|_{C^{\gamma}(\overline{D(t)})}|h|^{\gamma}\leq C(\|u_0\|_{H^{\frac12+\ep}},T)\|\theta_0\|_{C^\gamma(\overline{D}_0)}|h|^{\gamma},
$$
by using Theorem \ref{Case2}. The term $J_1$ has to be decompose in the following manner $J_1=L_1+L_2+L_3+L_4+L_5+L_6$, where

$$
L_1=\int_{\{|x-y|<2|h|\}\cap D(t)}k_4(x+h-y)(\theta(x+h,t)-\theta(y,t))dy,
$$
$$L_2=-\int_{\{|x-y|<2|h|\}\cap D(t)}k_4(x-y)(\theta(x,t)-\theta(y,t))dy,
$$
$$L_3=\int_{\{|x-y|\geq 2|h|\}\cap D(t)}k_4(x-y)(\theta(x+h,t)-\theta(x,t))dy,
$$
$$L_4=\int_{\{|x-y|\geq 2|h|\}\cap D(t)}(k_4(x+h-y)-k_4(x-y))(\theta(x+h,t)-\theta(y,t))dy,
$$
$$L_5=-(\theta(x+h,t)-\theta(x,t))\pv\int_{D(t)}k_4(x-y)dy,
$$
$$L_6=-\theta(x+h,t)\int_{D(t)}(k_4(x+h-y)-k_4(x-y))dy.
$$
Then
\begin{equation*}
|L_1|\leq  \|\theta(t)\|_{C^\gamma(\overline{D(t)})}\int_{|x-y|<2|h|} |k_4(x+h-y)||x+h-y|^\gamma dy\leq c(\|u_0\|_{H^{\frac12+\varepsilon}},T)\|\theta_0\|_{C^\gamma(\overline{D}_0)}|h|^\gamma,
\end{equation*}
and analogously
\begin{equation*}
|L_2|\leq  c(\|u_0\|_{H^{\frac12+\varepsilon}},T)\|\theta_0\|_{C^\gamma(\overline{D}_0)}|h|^\gamma.
\end{equation*}
Adding $L_3$ and $L_5$ we find
\begin{equation*}
\begin{aligned}
|L_3+L_5|&\leq |\theta(x+h,t)-\theta(x,t)|\Big| \pv\int_{\{|x-y|<2|h|\}\cap D(t)}k_4(x-y)dy\Big|\\
&\leq c(\|u_0\|_{H^{\frac12+\varepsilon}},\delta,T)\|\theta_0\|_{C^\gamma(\overline{D}_0)}|h|^\gamma,
\end{aligned}	
\end{equation*}
where the principal value is bounded as in \eqref{L2} taking $2|h|$ smaller than the cutoff \eqref{cutoff}. The next bound is performed as follows
\begin{equation*}
\begin{aligned}
|L_4|&\leq \|\theta(t)\|_{C^\gamma(\overline{D(t)})}\int_{\{|x-y|\geq 2|h|\}\cap D(t)}|k_4(x+h-y)-k_4(x-y)||x+h-y|^\gamma dy\\
&\leq c(\|u_0\|_{H^{\frac12+\varepsilon}},T)\|\theta_0\|_{C^\gamma(\overline{D}_0)} \int_{|x-y|\geq 2|h|}|k_4(x+h-y)-k_4(x-y)||x-y|^\gamma dy\\
&\leq c(\|u_0\|_{H^{\frac12+\varepsilon}},T)\|\theta_0\|_{C^\gamma(\overline{D}_0)} \int_{|x-y|\geq 2|h|} \frac{|h|}{|x-y|^{4-\gamma}}dy\leq c(\|u_0\|_{H^{\frac12+\varepsilon}},T)\|\theta_0\|_{C^\gamma(\overline{D}_0)} |h|^\gamma,
\end{aligned}
\end{equation*}
so that it remains to bound $L_6$ to be done with the H\"older regularity for $I_1$, and consequently with the proof. Therefore, the rest of the section is dedicated to bound $L_6$.
\vspace{0.2cm}

 \noindent \underline{\textbf{Bounding $L_6$}}:
 \vspace{0.2cm}
 
 \noindent The goal is to take advantage of the fact that the kernels $k_4$ can be written as derivatives, so that one can integrate by parts to obtain operators on the boundary. Thus, we rewrite the term to see that
\begin{equation*}
\begin{aligned}
L_6&=-\theta(x+h,t)(k_4*1_{D(t)}(x+h)-k_4*1_{D(t)}(x))\\
&=-\theta(x+h,t)((k_4+\frac15\delta_0)*1_{D(t)}(x+h)-(k_4+\frac15\delta_0)*1_{D(t)}(x)).
\end{aligned}
\end{equation*}
From \eqref{kernels} and \eqref{Ki}, we recall that 
\begin{equation*}
\mathcal{F}\Big(k_4+\frac15\delta_0\Big)(\xi)=\frac{\xi_1^4}{|\xi|^4}=\xi_1\frac{\xi_1^3}{|\xi|^4}, 
\end{equation*}
so that
\begin{equation}\label{Gamma4}
k_4+\frac15\delta_0=\partial_1 \mathcal{F}^{-1}\Big(-i\frac{\xi_1^3}{|\xi|^4}\Big)=\partial_1 \left(-3\frac{x_1(x_2^2+x_3^2)}{|x|^5}\right)=\partial_1 \Gamma_4(x),
\end{equation}
in the distributional sense. Then, integration by parts gives that
\begin{equation*}
\begin{aligned}
L_6=-\theta(x+h,t)\int_{\partial D(t)} (\Gamma_4(x+h-y)-\Gamma_4(x-y))n_1(y,t)dS(y).
\end{aligned}
\end{equation*}
Now, all the terms that will appear are singular integrals in two variables with kernels that depend on the free surface.
Thus, we will use contour dynamics techniques together with the regularity previously obtained for the surface. Further splitting is needed to handle all singular integral kernels along the free surface.

For simplicity, from now on we disregard the dependence in time of the notation. 
We take a cutoff distance $\eta>0$ and denote $B_\eta=\{y\in \partial D(t): |x-y|<\eta\}$. We can always choose the size of $h$ small enough so that $x+h\in B_{\eta/2}$. Then we can write $L_6$ as follows 
\begin{equation}\label{L6}
\begin{aligned}
L_6&=-\theta(x+h)\int_{B_\eta} (\Gamma_4(x+h-y)-\Gamma_4(x-y))n_1(y)dS(y)\\
&\quad-\theta(x+h)\int_{\partial D(t)\smallsetminus B_\eta} (\Gamma_4(x+h-y)-\Gamma_4(x-y))n_1(y)dS(y)=M_1+M_2.
\end{aligned}
\end{equation}
For the $M_2$ term we use the mean value theorem to obtain
\begin{equation}\label{M2bound}
|M_2|\leq \frac{c}{\eta^3}\|\theta_0\|_{L^\infty}|\partial D(t)||h|.
\end{equation}
Therefore, it only remains to deal with the $M_1$ term. To estimate it we choose a parametrization $Z(\beta,t)$ on $B_\eta$ of the surface $\partial D(t)$ near the points $Z(\alpha,t)$, $Z(\tilde{\alpha},t)$ so that
\begin{equation*}
\begin{aligned}
M_1&=-\theta(x+h)\int_{Z^{-1}(B_\eta)} (\Gamma_4(Z(\alpha)-Z(\beta))-\Gamma_4(Z(\tilde{\alpha})-Z(\beta)))N_1(Z(\beta))d\beta,
\end{aligned}
\end{equation*}
and substituting $\Gamma_4$ from \eqref{Gamma4} we can write
\begin{equation*}
\begin{aligned}
M_1&\!=\!-\theta(x\!+\!h)\!\int_{Z^{-1}(B_\eta)}\!\!\!\!\!\!\!\!\!\frac{(Z_2(\alpha)\!-\!Z_2(\beta))^2\!+\!(Z_3(\alpha)\!-\!Z_3(\beta))^2}{|\alpha-\beta|^2}\frac{|\alpha-\beta|^5N_1(Z(\beta))}{|Z(\alpha)\!-\!Z(\beta)|^5}\frac{Z_1(\alpha)\!-\!Z_1(\beta)}{|\alpha-\beta|^3}d\beta\\
&\quad+\theta(x\!+\!h)\!\int_{Z^{-1}(B_\eta)}\!\!\!\!\!\!\!\!\!\frac{(Z_2(\tilde{\alpha})\!-\!Z_2(\beta))^2\!+\!(Z_3(\tilde{\alpha})\!-\!Z_3(\beta))^2}{|\tilde{\alpha}-\beta|^2}\frac{|\tilde{\alpha}-\beta|^5N_1(Z(\beta))}{|Z(\tilde{\alpha})\!-\!Z(\beta)|^5}\frac{Z_1(\tilde{\alpha})\!-\!Z_1(\beta)}{|\tilde{\alpha}-\beta|^3}d\beta.
\end{aligned}
\end{equation*}
For convenience, we will choose the parametrization with isothermal coordinates (see e.g. \cite{Taylor11}, Chap. 5.10), i.e., verifying that
\begin{equation}\label{isothermal}
\partial_{\alpha_1}Z(\alpha,t)\cdot \partial_{\alpha_2}Z(\alpha,t)=0,\hspace{1cm}|\partial_{\alpha_1}Z(\alpha,t)|^2=|\partial_{\alpha_2}Z(\alpha,t)|^2.
\end{equation}
Then, we add and subtract the appropriate quantities and group the terms together,
\begin{equation}\label{M1}
M_1=O_1(\alpha)-O_1(\tilde{\alpha})+O_2(\alpha)-O_2(\tilde{\alpha}),
\end{equation}
where
\begin{equation*}
O_1(\alpha)=-\theta(x+h)\pv\!\int_{Z^{-1}(B_\eta)}\!\!\!\!\!\!\!\!\!\!\!\!\!\frac{(\partial_\alpha Z_2(\alpha)\!\cdot\!(\alpha\!-\!\beta))^2\!+\!(\partial_\alpha Z_3(\alpha)\!\cdot\!(\alpha\!-\!\beta))^2}{|\alpha-\beta|^2}\frac{N_1(Z(\alpha))}{|\partial_{\alpha_1}Z(\alpha)|^5}\frac{Z_1(\alpha)\!-\!Z_1(\beta)}{|\alpha-\beta|^3}d\beta,
\end{equation*}
\begin{equation*}
\begin{aligned}
O_2(\alpha)&=\!-\theta(x+h)\!\int_{Z^{-1}(B_\eta)}\!\!\!\!\!\!\!\!\!\!\!\!\frac{(Z_2(\alpha)\!-\!Z_2(\beta))^2\!+\!(Z_3(\alpha)\!-\!Z_3(\beta))^2}{|\alpha-\beta|^2}\frac{|\alpha-\beta|^5N_1(Z(\beta))}{|Z(\alpha)\!-\!Z(\beta)|^5}\frac{Z_1(\alpha)\!-\!Z_1(\beta)}{|\alpha-\beta|^3} d\beta\\
&\quad+\theta(x+h)\pv\!\int_{Z^{-1}(B_\eta)}\!\!\!\!\!\!\!\!\!\!\!\!\!\!\!\frac{(\partial_\alpha Z_2(\alpha)\!\cdot\!(\alpha\!-\!\beta))^2\!+\!(\partial_\alpha Z_3(\alpha)\!\cdot\!(\alpha\!-\!\beta))^2}{|\alpha-\beta|^2}\frac{N_1(Z(\alpha))}{|\partial_{\alpha_1} Z(\alpha)|^5}\frac{Z_1(\alpha)\!-\!Z_1(\beta)}{|\alpha-\beta|^3}d\beta.
\end{aligned}
\end{equation*}
To conclude the proof, we deal with the bounds for  $O_1(\alpha)-O_1(\alpha)$ and the bounds for  $O_2(\alpha)-O_2(\alpha)$ in two different subsections.

\vspace{0.3cm}
\noindent \underline{Bounding $O_1(\alpha)-O_1(\alpha)$}:
\vspace{0.3cm}

\noindent The term $O_1$ can be decomposed further:
\begin{equation}\label{O1}
O_1(\alpha)=\sum_{j=2,3} \left(P_{1,j}(\alpha)+P_{2,j}(\alpha)+P_{3,j}(\alpha)\right),
\end{equation}
where
\begin{equation*}
\begin{aligned}
P_{1,j}(\alpha)&=-\theta(x+h)\frac{(\partial_{\alpha_1}Z_j(\alpha))^2N_1(Z(\alpha))}{|\partial_{\alpha_1}Z(\alpha)|^5}\pv\int_{Z^{-1}(B_\eta)}\frac{(\alpha_1-\beta_1)^2}{|\alpha-\beta|^5}(Z_1(\alpha)\!-\!Z_1(\beta))d\beta,\\
P_{2,j}(\alpha)&=-\theta(x+h)\frac{(\partial_{\alpha_2}Z_j(\alpha))^2N_1(Z(\alpha))}{|\partial_{\alpha_1}Z(\alpha)|^5}\pv\int_{Z^{-1}(B_\eta)}\frac{(\alpha_2-\beta_2)^2}{|\alpha-\beta|^5}(Z_1(\alpha)\!-\!Z_1(\beta))d\beta,\\
P_{3,j}(\alpha)&=\!-2\theta(x\!+\!h)\frac{\partial_{\alpha_1}Z_j(\alpha)\partial_{\alpha_2}Z_j(\alpha)N_1(Z(\alpha))}{|\partial_{\alpha_1}Z(\alpha)|^5}\pv\!\!\int_{Z^{-1}(B_\eta)}\!\!\!\!\!\!\!\!\!\!\frac{(\alpha_1\!-\!\beta_1)(\alpha_2\!-\!\beta_2)}{|\alpha-\beta|^5}(Z_1(\alpha)\!-\!Z_1(\beta))d\beta.
\end{aligned}
\end{equation*}
We define the following quantity
\begin{equation*}
F(Z)(\alpha,\beta,t)=\frac{|\alpha-\beta|}{|Z(\alpha,t)-Z(\beta,t)|},  \hspace{0.5cm}\alpha,\beta\in Z^{-1}(B_\eta),
\end{equation*}
which measures the lack of self-intersection of $Z(t)$ on $B_\eta$. Then, it is not difficult to see that
\begin{equation}\label{P1}
|P_{1,j}(\alpha)-P_{1,j}(\tilde{\alpha})|\leq Q_1+Q_2,
\end{equation}
with
\begin{equation*}
Q_1=c\|\theta_0\|_{L^\infty}\|F(Z)\|_{L^\infty}^5\|\partial_\alpha Z\|_{L^\infty}^3 \|\partial_\alpha Z\|_{C^\gamma}|\alpha-\tilde{\alpha}|^\gamma \Big|\int_{Z^{-1}(B_\eta)}\frac{(\tilde{\alpha}_1-\beta_1)^2}{|\tilde{\alpha}-\beta|^5}(Z_1(\tilde{\alpha})-Z_1(\beta))d\beta\Big|,
\end{equation*}
\begin{multline*}
Q_2=c\|\theta_0\|_{L^\infty}\|\partial_\alpha Z\|_{L^\infty}^4 \|F(Z)\|_{L^\infty}^5 \Big|\!\int_{Z^{-1}(B_\eta)}\!\!\!\!\!\!\!\!\!\!\!\!\!\!\!\!\!\frac{(Z_1(\alpha)\!-\!Z_1(\beta))(\alpha_1\!-\!\beta_1)^2}{|\alpha-\beta|^5}d\beta\!\\
-\!\!\int_{Z^{-1}(B_\eta)}\!\!\!\!\!\!\!\!\!\!\!\!\!\!\!\!\frac{(Z_1(\tilde{\alpha})\!-\!Z_1(\beta))(\tilde{\alpha}_1\!-\!\beta_1)^2}{|\tilde{\alpha}-\beta|^5}d\beta \Big|.
\end{multline*}
To deal with $Q_1$, we first notice that
\begin{equation*}
\begin{aligned}
\Big|\!\int_{Z^{-1}(B_\eta)}\!\!\!\!\!\!\!\!\!\!\!\!\!\!\frac{(Z_1(\tilde{\alpha})\!-\!Z_1(\beta))(\tilde{\alpha}_1-\beta_1)^2}{|\tilde{\alpha}-\beta|^5}d\beta\Big|&\!\leq\! \Big|\!\int_{Z^{-1}(B_\eta)}\!\!\!\!\!\!\!\!\!\!\!\!\!\!\!\!\frac{(\tilde{\alpha}_1\!-\!\beta_1)^2(Z_1(\tilde{\alpha})-Z_1(\beta)-\partial_{\alpha} Z_1(\tilde{\alpha})\cdot (\tilde{\alpha}-\beta))}{|\tilde{\alpha}-\beta|^5}d\beta\Big|\\
&\quad+\sum_{i=1,2} \Big|\int_{Z^{-1}(B_\eta)}\!\!\!\!\!\!\frac{(\tilde{\alpha}_1-\beta_1)^2(\tilde{\alpha}_i-\beta_i)}{|\tilde{\alpha}-\beta|^5}\partial_{\alpha_i} Z_1(\tilde{\alpha})d\beta\Big|,
\end{aligned}
\end{equation*}
so using that $Z(t)\in C^{1+\mu}$, $0<\mu<1$, we obtain
\begin{equation*}
\begin{aligned}
\Big|\!\int_{Z^{-1}(B_\eta)}\!\!\!\!\!\!\!\!\!\frac{(Z_1(\tilde{\alpha})\!-\!Z_1(\beta))(\tilde{\alpha}_1\!-\!\beta_1)^2}{|\tilde{\alpha}-\beta|^5}d\beta\Big|&\leq c\|Z\|_{C^{1+\mu}}\!+\! \sum_{i=1,2}\Big|\partial_{\alpha_i} Z_1(\tilde{\alpha})\int_{Z^{-1}(B_\eta)}\!\!\!\!\!\!\!\!\!\!\frac{(\tilde{\alpha}_1\!-\!\beta_1)^2(\tilde{\alpha}_i\!-\!\beta_i)}{|\tilde{\alpha}-\beta|^5}d\beta\Big|.
\end{aligned}
\end{equation*}
Recalling that $\tilde{\alpha}\in Z^{-1}(B_{\eta/2})$, we have that $d(\tilde{\alpha},Z^{-1}(\partial B_\eta))\geq \eta /(2\|\partial_{\alpha}Z\|_{L^\infty})$. Therefore, if we take $\delta=\eta/(4\|\partial_{\alpha}Z\|_{L^\infty})$, we find that
\begin{multline*}
\Big|\int_{Z^{-1}(B_\eta)}\frac{(\tilde{\alpha}_1-\beta_1)^2(\tilde{\alpha}_i-\beta_i)}{|\tilde{\alpha}-\beta|^5}d\beta\Big|=\Big|\underbrace{\int_{|\tilde{\alpha}-\beta|\leq \delta}\frac{(\tilde{\alpha}_1-\beta_1)^2(\tilde{\alpha}_i-\beta_i)}{|\tilde{\alpha}-\beta|^5}d\beta}_{=0}\\
\quad+\int_{Z^{-1}(B_\eta)\smallsetminus\{|\tilde{\alpha}-\beta|\leq \delta\}}\frac{(\tilde{\alpha}_1-\beta_1)^2(\tilde{\alpha}_i-\beta_i)}{|\tilde{\alpha}-\beta|^5}d\beta\Big|\leq c\frac{\|\partial_{\alpha}Z\|_{L^\infty}}{\eta},
\end{multline*}
and thus
\begin{equation*}
\begin{aligned}
\Big|\!\int_{Z^{-1}(B_\eta)}\frac{(Z_1(\tilde{\alpha})-Z_1(\beta))(\tilde{\alpha}_1-\beta_1)^2}{|\tilde{\alpha}-\beta|^5}d\beta\Big|&\leq c\left(\|Z\|_{C^{1+\mu}}+\frac{\|\partial_{\alpha}Z\|^2_{L^\infty}}{\eta}\right).
\end{aligned}
\end{equation*}
We conclude then the bound for $Q_1$
\begin{equation*}
\begin{aligned}
Q_1 &\leq c\|\theta_0\|_{L^\infty}\|F(Z)\|_{L^\infty}^5\|\partial_\alpha Z\|_{L^\infty}^3 \|\partial_\alpha Z\|_{C^\gamma}|\alpha-\tilde{\alpha}|^\gamma \left(\|Z\|_{C^{1+\mu}}+\frac{\|\partial_{\alpha}Z\|^2_{L^\infty}}{\eta}\right)\\
&\leq c\|\theta_0\|_{L^\infty}\|F(Z)\|_{L^\infty}^{5+\gamma}\|\partial_\alpha Z\|_{L^\infty}^3 \|\partial_\alpha Z\|_{C^\gamma} \left(\|Z\|_{C^{1+\mu}}+\frac{\|\partial_{\alpha}Z\|^2_{L^\infty}}{\eta}\right)|h|^\gamma,
\end{aligned}
\end{equation*}
that is, 
\begin{equation}\label{Q1bound}
Q_1\leq c(\|\theta_0\|_{L^\infty},\|F(Z)\|_{L^\infty},\|Z\|_{C^{1+\gamma}},\eta)|h|^\gamma.
\end{equation}
To deal with $Q_2$, we use the following identity
\begin{equation*}
-\frac{x_1^2}{|x|^5}=\partial_{x_1}\left(\frac{x_1^3}{|x|^5}\right)+\partial_{x_2}\left(\frac{x_1^2x_2}{|x|^5}\right), \hspace{0.2cm}x\in \mathbb{R}^2,
\end{equation*}
followed by integration by parts to show that
\begin{equation*}
\begin{aligned}
\int_{Z^{-1}(B_\eta)}\!\!\!\!\!\!\!\!\!\!\!\!\!\frac{(Z_1(\alpha)-Z_1(\beta))(\alpha_1-\beta_1)^2}{|\alpha-\beta|^5}d\beta&=\int_{Z^{-1}(B_\eta)}\!\!\!\!\!\!\!\!\!\!\frac{(\alpha_1\!-\!\beta_1)^3\partial_{\alpha_1} Z_1(\beta)}{|\alpha-\beta|^5}d\beta\\
&\quad+\int_{Z^{-1}(B_\eta)}\!\!\!\!\!\!\!\!\!\frac{(\alpha_1\!-\!\beta_1)^2(\alpha_2\!-\!\beta_2)}{|\alpha-\beta|^5}\partial_{\alpha_2} Z_1(\beta)d\beta\\
&\quad-\int_{\partial Z^{-1}(B_\eta)}\!\!\!\!\!\!\!\!\!\!\!\!\!\!\!\frac{(\alpha_1\!-\!\beta_1)^2(\alpha_2\!-\!\beta_2)(Z_1(\alpha)\!-\!Z_1(\beta))n_2(\beta)}{|\alpha-\beta|^5}dl(\beta)\\
&\quad-\int_{\partial Z^{-1}(B_\eta)}\!\!\!\!\!\!\!\!\!\!\!\!\frac{(\alpha_1\!-\!\beta_1)^3(Z_1(\alpha)\!-\!Z_1(\beta))n_1(\beta)}{|\alpha-\beta|^5}dl(\beta).
\end{aligned}
\end{equation*}
Then, $Q_2$ is bounded as follows
\begin{equation}\label{Q2}
Q_2\leq c\|\theta_0\|_{L^\infty}\|\partial_\alpha Z\|_{L^\infty}^4 \|F(Z)\|_{L^\infty}^5 (R_1+R_2+R_3+R_4),
\end{equation}
where
$$
R_1=\Big|\int_{Z^{-1}(B_\eta)}\frac{(\alpha_1-\beta_1)^3\partial_{\alpha_1}Z_1(\beta)}{|\alpha-\beta|^5}d\beta-\int_{Z^{-1}(B_\eta)}\frac{(\tilde{\alpha}_1-\beta_1)^3\partial_{\alpha_1} Z_1(\beta)}{|\tilde{\alpha}-\beta|^5}d\beta\Big|,
$$
$$
R_2=\Big|\!\int_{\partial Z^{-1}(B_\eta)}\!\!\!\!\!\!\!\!\!\!\!\!\!\!\!\frac{(\alpha_1-\beta_1)^3(Z_1(\alpha)\!-\!Z_1(\beta))n_1(\beta)}{|\alpha-\beta|^5}dl(\beta)-\!
\int_{\partial Z^{-1}(B_\eta)}\!\!\!\!\!\!\!\!\!\!\!\!\!\!\!\frac{(\tilde{\alpha}-\beta_1)^3(Z_1(\tilde{\alpha})\!-\!Z_1(\beta))n_1(\beta)}{|\tilde{\alpha}-\beta|^5}dl(\beta)\Big|,
$$
$$
R_3=\Big|\int_{Z^{-1}(B_\eta)}\!\!\!\!\!\!\!\!\frac{(\alpha_1-\beta_1)^2(\alpha_2-\beta_2)}{|\alpha-\beta|^5}\partial_{\alpha_1} Z_1(\beta)d\beta-\int_{Z^{-1}(B_\eta)}\!\!\!\!\!\!\!\!\frac{(\tilde{\alpha}_1-\beta_1)^2(\tilde{\alpha}_2-\beta_2)}{|\tilde{\alpha}-\beta|^5}\partial_{\alpha_1} Z_1(\beta)d\beta\Big|,
$$
and
\begin{equation*}
\begin{aligned}
R_4=\Big|\!\int_{\partial Z^{-1}(B_\eta)}&\frac{(\alpha_1-\beta_1)^2(\alpha_2-\beta_2)(Z_1(\alpha)\!-\!Z_1(\beta))n_2(\beta)}{|\alpha-\beta|^5}dl(\beta)\\
&-\int_{\partial Z^{-1}(B_\eta)}\frac{(\tilde{\alpha}_1-\beta_1)^2(\tilde{\alpha}_2-\beta_2)(Z_1(\tilde{\alpha})\!-\!Z_1(\beta))n_2(\beta)}{|\tilde{\alpha}-\beta|^5}dl(\beta)\Big|.
\end{aligned}
\end{equation*}
Introducing one more splitting, the term $R_1$ is written in the following manner
\begin{equation*}
\begin{aligned}
R_1=\big|S_1+S_2+S_3+S_4+S_5+S_6\big|,
\end{aligned}
\end{equation*}
where
$$
S_1=\int_{Z^{-1}(B_\eta)\cap\{|\alpha-\beta|<2|\alpha-\tilde{\alpha}|\}}\frac{(\alpha_1-\beta_1)^3}{|\alpha-\beta|^5}(\partial_{\alpha_1} Z_1(\alpha)-\partial_{\alpha_1}Z_1(\beta))d\beta,
$$
$$S_2=-\int_{Z^{-1}(B_\eta)\cap\{|\tilde{\alpha}-\beta|\leq 2|\alpha-\tilde{\alpha}|\}}\frac{(\tilde{\alpha}_1-\beta_1)^3}{|\tilde{\alpha}-\beta|^5}(\partial_{\alpha_1}Z_1(\tilde{\alpha})-\partial_{\alpha_1}Z(\beta))d\beta,$$
$$S_3=\int_{Z^{-1}(B_\eta)\cap\{|\alpha-\beta|\geq 2|\alpha-\tilde{\alpha}|\}}\frac{(\alpha_1-\beta_1)^3}{|\alpha-\beta|^5}(\partial_{\alpha_1}Z_1(\alpha)-\partial_{\alpha_1}Z_1(\tilde{\alpha}))d\beta,$$
$$S_4=\int_{Z^{-1}(B_\eta)\cap\{|\tilde{\alpha}-\beta|\geq 2|\alpha-\tilde{\alpha}|\}}\left(\frac{(\alpha_1-\beta_1)^3}{|\alpha-\beta|^5}-\frac{(\tilde{\alpha}_1-\beta_1)^3}{|\tilde{\alpha}-\beta|^5}\right)(\partial_{\alpha_1}Z_1(\tilde{\alpha})-\partial_{\alpha_1}Z(\beta))d\beta,$$
$$S_5=-\int_{Z^{-1}(B_\eta)}\frac{(\alpha_1-\beta_1)^3}{|\alpha-\beta|^5}(\partial_{\alpha_1}Z_1(\alpha)-\partial_{\alpha_1}Z_1(\tilde{\alpha}))d\beta,$$
$$S_6=\int_{Z^{-1}(B_\eta)}\left(\frac{(\tilde{\alpha}_1-\beta_1)^3}{|\tilde{\alpha}-\beta|^5}-\frac{(\alpha_1-\beta_1)^3}{|\alpha-\beta|^5}\right)\partial_{\alpha_1}Z_1(\tilde{\alpha})d\beta.
$$
One immediately obtains that
\begin{equation*}
|S_1|+|S_2|+|S_4|\leq c \|Z\|_{C^{1+\gamma}}|\alpha-\tilde{\alpha}|^\gamma\leq c \|Z\|_{C^{1+\gamma}}\|F(Z)\|_{L^\infty}^\gamma |h|^\gamma. 
\end{equation*}
The terms $S_3$ and $S_5$ are joined together
\begin{equation*}
S_3+S_5=(\partial_{\alpha_1}Z_1(\alpha)-\partial_{\alpha_1}Z_1(\tilde{\alpha}))\int_{Z^{-1}(B_\eta)\cap\{|\alpha-\beta|\leq 2|\alpha-\tilde{\alpha}|\}}\frac{(\alpha_1-\beta_1)^3}{|\alpha-\beta|^5}d\beta.
\end{equation*}
Recalling that $Z(\alpha)$ is the center of $B_\eta$, we know that $d(\alpha,\partial Z^{-1}(B_\eta))\geq \eta/\|\partial_\alpha Z\|_{L^\infty}$. Then, since $|\alpha-\tilde{\alpha}|\leq \|F(Z)\|_{L^\infty}|h|$, we can choose $|h|<(\eta/2)/(\|\partial_{\alpha}Z\|_{L^\infty}\|F(Z)\|_{L^\infty})$ to guarantee that
\begin{equation*}
|\alpha-\tilde{\alpha}|< \frac{\eta}{2\|\partial_{\alpha} Z\|_{L^\infty}},
\end{equation*}
so that the integral is on a disk and therefore vanishes.
Finally, integration by parts in $S_6$ shows that 
\begin{equation*}
S_6\!=\!\partial_{\alpha_1}Z_1(\tilde{\alpha}) \int_{\partial Z^{-1}(B_\eta)}\!\!\left(\frac{3(\tilde{\alpha}_1\!-\!\beta_1)^2\!+\!2(\tilde{\alpha}_2\!-\!\beta_2)^2}{3|\tilde{\alpha}-\beta|^3}\!-\!\frac{3(\alpha_1\!-\!\beta_1)^2\!+\!2(\alpha_2\!-\!\beta_2)^2}{3|\alpha-\beta|^3}\right)n_1(\beta)dl(\beta),
\end{equation*}
so, since $\alpha,\tilde{\alpha}\in Z^{-1}(B_{\eta/2})$, we can apply the mean value theorem to conclude that
\begin{equation*}
|S_{6}|\leq c\|\partial_{\alpha} Z\|_{L^\infty}\left(\frac{2\|\partial_{\alpha}Z\|_{L^\infty}}{\eta}\right)^2|\alpha-\tilde{\alpha}|\leq c \frac{\|\partial_{\alpha}Z\|_{L^\infty}^3}{\eta^2}\|F(Z)\|_{L^\infty}|h|.
\end{equation*}
Joining the above bounds we have that
\begin{equation}\label{R1bound}
R_1\leq c\left(\|Z\|_{C^{1+\gamma}}\|F(Z)\|_{L^\infty}^\gamma+ \frac{\|\partial_{\alpha}Z\|_{L^\infty}^3}{\eta^2}\|F(Z)\|_{L^\infty}\right)|h|^\gamma.
\end{equation}
We now rewrite $R_2$ as follows
\begin{multline*}
R_2=\Big|\int_{\partial Z^{-1}(B_\eta)}\frac{(\alpha_1-\beta_1)^3}{|\alpha-\beta|^5}(Z_1(\alpha)-Z_1(\tilde{\alpha}))n_1(\beta)dl(\beta)\\
-\int_{\partial Z^{-1}(B_\eta)}\left(\frac{(\tilde{\alpha}_1-\beta_1)^3}{|\tilde{\alpha}-\beta|^5}-\frac{(\alpha_1-\beta_1)^3}{|\alpha-\beta|^5}\right)(Z_1(\tilde{\alpha})-Z_1(\beta))n_1(\beta)dl(\beta)\Big|,
\end{multline*}
which recalling again that $\alpha,\tilde{\alpha}\in Z^{-1}(B_{\eta/2})$ can be bounded by the following
\begin{equation}\label{R2bound}
R_2\leq c \frac{\|\partial_\alpha Z\|_{L^\infty}^2}{\eta^2}\|Z\|_{C^\gamma}\|F(Z)\|_{L^\infty}^\gamma |h|^\gamma + c\|Z\|_{L^\infty}\frac{\|\partial_{\alpha}Z\|_{L^\infty}^3}{\eta^3}\|F(Z)\|_{L^\infty}|h|.
\end{equation}
The terms $R_3$ and $R_4$ can be bounded analogously to $R_1$ and $R_2$, respectively. Introducing the bounds \eqref{R1bound}, \eqref{R2bound} back in \eqref{Q2}, we obtain that
\begin{equation}\label{Q2bound}
Q_2\leq c(\|\theta_0\|_{L^\infty},\|F(Z)\|_{L^\infty},\|Z\|_{C^{1+\gamma}},\eta)|h|^\gamma.
\end{equation}
From \eqref{P1}, the bounds \eqref{Q1bound} and \eqref{Q2bound}  yields that
\begin{equation}\label{P1bound}
|P_{1,j}(\alpha)-P_{1,j}(\tilde{\alpha})|\leq c(\|\theta_0\|_{L^\infty},\|F(Z)\|_{L^\infty},\|Z\|_{C^{1+\gamma}},\eta)|h|^\gamma.
\end{equation} 
The terms $P_{2,j}$ and $P_{3,j}$ can be estimated analogously, so we conclude in \eqref{O1} that
\begin{equation}\label{O1bound}
|O_1(\alpha)-O_1(\tilde{\alpha})|\leq c(\|\theta_0\|_{L^\infty},\|F(Z)\|_{L^\infty},\|Z\|_{C^{1+\gamma}},\eta)|h|^\gamma.
\end{equation} 
Thus, it only remains to bound the term $O_2(\alpha)-O_2(\alpha)$.

\vspace{0.3cm}
\noindent \underline{Bounding $O_2(\alpha)-O_2(\alpha)$}:
\vspace{0.3cm}

\noindent Going back to \eqref{M1}, for simplicity of notation we will denote
\begin{equation*}
G(\alpha,\beta)=\frac{(Z_2(\alpha)-Z_2(\beta))^2+(Z_3(\alpha)-Z_3(\beta))^2}{|\alpha-\beta|^2}\frac{|\alpha-\beta|^5N_1(Z(\beta))}{|Z(\alpha)-Z(\beta)|^5},
\end{equation*}
and define
\begin{equation*}
\mathcal{G}(\alpha,\beta)=\frac{(\partial_\alpha Z_2(\alpha)\!\cdot\!(\alpha\!-\!\beta))^2\!+\!(\partial_\alpha Z_3(\alpha)\!\cdot\!(\alpha\!-\!\beta))^2}{|\alpha-\beta|^2}\frac{N_1(Z(\alpha))}{|\partial_{\alpha}Z(\alpha)|^5}.
\end{equation*}
Then, by using the mean value theorem, the term $O_2(\alpha)$ is rewritten as follows
\begin{equation*}
O_2(\alpha)\!=\!-\theta(x+h)\!\int_0^1\!\!\int_{Z^{-1}(B_\eta)}\! \frac{\alpha\!-\!\beta}{|\alpha\!-\!\beta|^3}\cdot(\partial_{\alpha}Z_1((1\!-\!r)\beta\!+\!r\alpha)) (G(\alpha,\beta)\!-\!\mathcal{G}(\alpha,\beta)) d\beta dr.
\end{equation*}
Therefore, the difference $|O_2(\alpha)-O_2(\tilde{\alpha})|$ can be bounded by introducing the following splitting
\begin{equation}\label{O2}
\begin{aligned}
|O_2(\alpha)-O_2(\tilde{\alpha})|\leq \|\theta_0\|_{L^\infty} (P_4+P_5+P_6+P_7+P_8),
\end{aligned}
\end{equation}
where
\begin{equation*}
P_4=\Big|\!\int_0^1\!\!\int_{Z^{-1}(B_\eta)\cap \{|\alpha-\beta|\leq 2|\alpha-\tilde{\alpha}|\}}\! \frac{\alpha\!-\!\beta}{|\alpha\!-\!\beta|^3}\cdot\partial_{\alpha}Z_1((1\!-\!r)\beta\!+\!r\alpha) (G(\alpha,\beta)\!-\!\mathcal{G}(\alpha,\beta)) d\beta dr \Big|,
\end{equation*}
\begin{equation*}
P_5=\Big|\!\int_0^1\!\!\int_{Z^{-1}(B_\eta)\cap \{|\alpha-\beta|\leq 2|\alpha-\tilde{\alpha}|\}}\! \frac{\tilde{\alpha}\!-\!\beta}{|\tilde{\alpha}\!-\!\beta|^3}\cdot\partial_{\alpha}Z_1((1\!-\!r)\beta\!+\!r\tilde{\alpha}) (G(\tilde{\alpha},\beta)\!-\!\mathcal{G}(\tilde{\alpha},\beta)) d\beta dr \Big|,
\end{equation*}
\begin{equation*}
P_6\!=\!\Big|\!\int_0^1\!\!\int_{Z^{-1}(B_\eta)\cap \{|\alpha-\beta|\geq2|\alpha-\tilde{\alpha}|\}}\!\!\!\!\!\!\!\!\!\!\!\!\!\!\!\!\!\!\!\!\!\!\!\!\!\!\!\!\!\!\!\!\!\!\!\!\!\!\!\!\!\!\!\!\!\!\!\!\!\!\!\!\! \frac{\big(\partial_{\alpha}Z_1((1\!-\!r)\beta\!+\!r\alpha)-\partial_{\alpha}Z_1((1\!-\!r)\beta\!+\!r\tilde{\alpha})\big)\cdot(\tilde{\alpha}-\beta) (G(\tilde{\alpha},\beta)-\mathcal{G}(\tilde{\alpha},\beta))}{|\tilde{\alpha}\!-\!\beta|^3} d\beta dr \Big|,
\end{equation*}
\begin{equation*}
P_7=\Big|\!\int_0^1\!\!\int_{Z^{-1}(B_\eta)\cap \{|\alpha-\beta|\geq 2|\alpha-\tilde{\alpha}|\}}\! \!\!\!\!\!\!\!\!\!\!\!\!\!\!\!\!\!\!\!\!\!\!\!\!\!\!\!\!\!\!\!\!\!\!\!\partial_{\alpha}Z_1((1-r)\beta+r\alpha)\cdot\left(\frac{\alpha\!-\!\beta}{|\alpha\!-\!\beta|^3}\!-\!\frac{\tilde{\alpha}\!-\!\beta}{|\tilde{\alpha}\!-\!\beta|^3}\right)  (G(\tilde{\alpha},\beta)-\mathcal{G}(\tilde{\alpha},\beta)) d\beta dr \Big|,
\end{equation*}
\begin{equation*}
P_8\!=\!\Big|\!\int_0^1\!\!\!\int_{Z^{-1}(B_\eta)\cap \{|\alpha-\beta|\geq 2|\alpha-\tilde{\alpha}|\}}\! \!\!\!\!\!\!\!\!\!\!\!\!\!\!\!\!\!\!\!\!\!\!\!\!\!\!\!\!\!\!\!\!\!\!\!\!\!\!\!\!\!\!\!\frac{(\alpha\!-\!\beta)\cdot\partial_{\alpha}Z_1((1-r)\beta+r\alpha)\big(G(\alpha,\beta)\!-\!G(\tilde{\alpha},\beta)\!-\!\mathcal{G}(\alpha,\beta)\!+\!\mathcal{G}(\tilde{\alpha},\beta)\big) }{|\alpha\!-\!\beta|^3}  d\beta dr \Big|.
\end{equation*}
To deal with the first two terms, we first notice that for both $w=\alpha$ and $w=\tilde{\alpha}$ it holds that
\begin{equation*}
\begin{aligned}
|G(w,\beta)-\mathcal{G}(w,\beta)|&\leq c\left(1+ \|\partial_{\alpha}Z\|_{L^\infty}\|F(Z)\|_{L^\infty}\right)\|\partial_{\alpha}Z\|_{L^\infty}^3\|F(Z)\|_{L^\infty}^5\|\partial_{\alpha}Z\|_{C^\gamma}|w-\beta|^\gamma\\
&\leq c(\|F(Z)\|_{L^\infty},\|Z\|_{C^{1+\gamma}})|w-\beta|^\gamma.
\end{aligned}
\end{equation*}
Therefore, one can integrate to find that 
\begin{equation}\label{P1P2}
P_4+P_5\leq c(\|F(Z)\|_{L^\infty},\|Z\|_{C^{1+\gamma}})|\alpha-\tilde{\alpha}|^\gamma\leq c(\|F(Z)\|_{L^\infty},\|Z\|_{C^{1+\gamma}})|h|^\gamma.
\end{equation}
Next, $P_6$ is readily bounded as follows
\begin{equation}\label{P3}
P_6\leq c(\|F(Z)\|_{L^\infty},\|Z\|_{C^{1+\gamma}})|h|^\gamma.
\end{equation}
The mean value theorem applied in $P_7$ provides that
\begin{equation*}
P_7\leq c(\|F(Z)\|_{L^\infty},\|Z\|_{C^{1+\gamma}})\|\partial_{\alpha}Z|\|_{L^\infty}\int_{Z^{-1}(B_\eta)\cap \{|\alpha-\beta|\geq 2|\alpha-\tilde{\alpha}|\}}\frac{|\alpha-\tilde{\alpha}|}{|\alpha-\beta|^3}|\tilde{\alpha}-\beta|^\gamma  d\beta, 
\end{equation*}
and since $|\tilde{\alpha}-\beta|\leq \frac32|\alpha-\beta|$, we conclude that
\begin{equation}\label{P4}
P_7\leq c(\|F(Z)\|_{L^\infty},\|Z\|_{C^{1+\gamma}})|h|^\gamma.
\end{equation}
It remains to deal with $P_8$. To bound it we decompose further $G(\alpha,\beta)-G(\tilde{\alpha},\beta)-\mathcal{G}(\alpha,\beta)+\mathcal{G}(\tilde{\alpha},\beta)$. First,
\begin{equation}\label{G1G2}
G(\alpha,\beta)=G_1(\alpha,\beta)+G_2(\alpha,\beta),
\end{equation}
with 
\begin{equation*}
\begin{aligned}
G_1(\alpha,\beta)=\frac{(Z_2(\alpha)\!-\!Z_2(\beta))^2}{|\alpha-\beta|^2}\frac{|\alpha\!-\!\beta|^5N_1(Z(\beta))}{|Z(\alpha)-Z(\beta)|^5},G_2(\alpha,\beta)=\frac{(Z_3(\alpha)\!-\!Z_3(\beta))^2}{|\alpha-\beta|^2}\frac{|\alpha\!-\!\beta|^5N_1(Z(\beta))}{|Z(\alpha)-Z(\beta)|^5}.
\end{aligned}
\end{equation*}
Analogously, 
\begin{equation}\label{G1G22}
\mathcal{G}(\alpha,\beta)=\mathcal{G}_1(\alpha,\beta)+\mathcal{G}_2(\alpha,\beta),
\end{equation}
where 
\begin{equation*}
\begin{aligned}
\mathcal{G}_1(\alpha,\beta)=\frac{(\partial_\alpha Z_2(\alpha)\!\cdot\!(\alpha\!-\!\beta))^2}{|\alpha-\beta|^2}\frac{N_1(Z(\alpha))}{|\partial_{\alpha}Z(\alpha)|^5},\hspace{0.5cm}
\mathcal{G}_2(\alpha,\beta)=\frac{(\partial_\alpha Z_3(\alpha)\!\cdot\!(\alpha\!-\!\beta))^2}{|\alpha-\beta|^2}\frac{N_1(Z(\alpha))}{|\partial_{\alpha}Z(\alpha)|^5}.
\end{aligned}
\end{equation*}
Then, 
\begin{equation}\label{G1H1H2}
G_1(\alpha,\beta)-G_1(\tilde{\alpha},\beta)=H_1+H_2,
\end{equation}
\begin{equation*}
\begin{aligned}
H_1&=\left(\frac{(Z_2(\alpha)\!-\!Z_2(\beta))^2}{|\alpha-\beta|^2}-\frac{(Z_2(\tilde{\alpha})\!-\!Z_2(\beta))^2}{|\tilde{\alpha}-\beta|^2}\right)\frac{|\alpha\!-\!\beta|^5N_1(Z(\beta))}{|Z(\alpha)-Z(\beta)|^5},\\
H_2&=\frac{(Z_2(\tilde{\alpha})\!-\!Z_2(\beta))^2}{|\tilde{\alpha}-\beta|^2}\left(\frac{|\alpha\!-\!\beta|^5}{|Z(\alpha)-Z(\beta)|^5}-\frac{|\tilde{\alpha}\!-\!\beta|^5}{|Z(\tilde{\alpha})-Z(\beta)|^5}\right)N_1(Z(\beta)).
\end{aligned}
\end{equation*}
Furthermore, 
\begin{equation*}
H_1=\left(\frac{Z_2(\alpha)\!-\!Z_2(\beta)}{|\alpha-\beta|}\!-\!\frac{Z_2(\tilde{\alpha})\!-\!Z_2(\beta)}{|\tilde{\alpha}-\beta|}\right)\left(\frac{Z_2(\alpha)\!-\!Z_2(\beta)}{|\alpha-\beta|}\!+\!\frac{Z_2(\tilde{\alpha})\!-\!Z_2(\beta)}{|\tilde{\alpha}-\beta|}\right)\frac{|\alpha\!-\!\beta|^5N_1(Z(\beta))}{|Z(\alpha)-Z(\beta)|^5}.
\end{equation*}
We can perform a similar decomposition of 
\begin{equation}\label{G1H1H22}
\mathcal{G}_1(\alpha,\beta)-\mathcal{G}_1(\tilde{\alpha},\beta)=\mathcal{H}_1+\mathcal{H}_2,
\end{equation}
where 
$$\mathcal{H}_1=\left(\frac{\partial_\alpha Z_2(\alpha)\!\cdot\!(\alpha\!-\!\beta)}{|\alpha-\beta|}\!-\!\frac{\partial_\alpha Z_2(\tilde{\alpha})\!\cdot\!(\tilde{\alpha}\!-\!\beta)}{|\tilde{\alpha}-\beta|}\right)\left(\frac{\partial_\alpha Z_2(\alpha)\!\cdot\!(\alpha\!-\!\beta)}{|\alpha-\beta|}\!+\!\frac{\partial_\alpha Z_2(\tilde{\alpha})\!\cdot\!(\tilde{\alpha}\!-\!\beta)}{|\tilde{\alpha}-\beta|}\right)\frac{N_1(Z(\alpha))}{|\partial_{\alpha}Z(\alpha)|^5},
$$
$$
\mathcal{H}_2=\frac{(\partial_\alpha Z_2(\tilde{\alpha})\!\cdot\!(\tilde{\alpha}\!-\!\beta))^2}{|\tilde{\alpha}-\beta|^2}\left(\frac{N_1(Z(\alpha))}{|\partial_{\alpha}Z(\alpha)|^5}-\frac{N_1(Z(\tilde{\alpha}))}{|\partial_{\alpha}Z(\tilde{\alpha})|^5}\right).
$$
Denote 
$$g(\alpha,\beta)=\frac{Z_2(\alpha)\!-\!Z_2(\beta)}{|\alpha-\beta|}, \hspace{0.5cm}\text{g}(\alpha,\beta)=\frac{\partial_\alpha Z_2(\alpha)\!\cdot\!(\alpha\!-\!\beta)}{|\alpha-\beta|}.$$
Then, we find that
\begin{equation}\label{H1H1}
H_1-\mathcal{H}_1=Y_1+Y_2,
\end{equation}
\begin{equation*}
\begin{aligned}
Y_1&=\big((g(\alpha,\beta)-g(\tilde{\alpha},\beta))-(\text{g}(\alpha,\beta)-\text{g}(\tilde{\alpha},\beta))\big)\big(g(\alpha,\beta)+g(\tilde{\alpha},\beta)\big)\frac{|\alpha\!-\!\beta|^5N_1(Z(\beta))}{|Z(\alpha)-Z(\beta)|^5},\\
Y_2&=(\text{g}(\alpha,\beta)-\text{g}(\tilde{\alpha},\beta))\Big(\big(g(\alpha,\beta)+g(\tilde{\alpha},\beta)\big)\frac{|\alpha\!-\!\beta|^5N_1(Z(\beta))}{|Z(\alpha)-Z(\beta)|^5} - \big(\text{g}(\alpha,\beta)+\text{g}(\tilde{\alpha},\beta) \big)\frac{N_1(Z(\alpha))}{|\partial_{\alpha}Z(\alpha)|^5} \Big).
\end{aligned}
\end{equation*}
We can bound $Y_1$ as follows
\begin{equation*}
|Y_1|\leq c(\|F(Z)\|_{L^\infty},\|\partial_{\alpha}Z\|_{L^\infty}) \big|(g(\alpha,\beta)-g(\tilde{\alpha},\beta))-(\text{g}(\alpha,\beta)-\text{g}(\tilde{\alpha},\beta))\big|,
\end{equation*}
where
\begin{equation*}
\big|g(\alpha,\beta)-g(\tilde{\alpha},\beta))-(\text{g}(\alpha,\beta)-\text{g}(\tilde{\alpha},\beta)\big|=\big|\delta g_1+\delta g_2+\delta g_3  \big|,
\end{equation*}
and, for $\sigma\in(0,1-\gamma)$, 
\begin{equation*}
|\delta g_1|=\Big|\int_0^1\Big(\partial_{\alpha}Z_2((1-r)\beta+r\alpha)-\partial_{\alpha}Z_2((1-r)\beta+r\tilde{\alpha})\Big)\cdot\frac{\alpha-\beta}{|\alpha-\beta|}dr\Big|\leq c\|\partial_{\alpha}Z_2\|_{C^{\gamma+\sigma}}|\alpha-\tilde{\alpha}|^{\gamma+\sigma},
\end{equation*}
\begin{equation*}
|\delta g_2|=\Big|-(\partial_{\alpha}Z_2(\alpha)-\partial_{\alpha}Z_2(\tilde{\alpha}))\cdot\frac{\alpha-\beta}{|\alpha-\beta|}\Big|\leq c\|\partial_{\alpha}Z_2\|_{C^{\gamma+\sigma}}|\alpha-\tilde{\alpha}|^{\gamma+\sigma},
\end{equation*}
\begin{equation*}
|\delta g_3|\!=\!\Big|\int_0^1\! \Big(\partial_{\alpha}Z_2((1-r)\beta+r\tilde{\alpha})-\partial_{\alpha}Z_2(\tilde{\alpha})\Big)\cdot\Big(\frac{\alpha\!-\!\beta}{|\alpha\!-\!\beta|}- \frac{\tilde{\alpha}\!-\!\beta}{|\tilde{\alpha}\!-\!\beta|} \Big)dr\Big|\leq c\|\partial_{\alpha}Z_2\|_{C^{\gamma+\sigma}}\frac{|\alpha-\tilde{\alpha}|}{|\alpha-\beta|^{1-\gamma-\sigma}}.
\end{equation*}
In the last inequality above, we used that in $P_8$ we are integrating in the region $|\alpha-\beta|\geq 2|\alpha-\tilde{\alpha}|$.
Thus, we have found the following bound for $Y_1$ 
\begin{equation*}
|Y_1|\leq c(\|F(Z)\|_{L^\infty},\|\partial_{\alpha}Z\|_{C^{\gamma+\sigma}})|\alpha-\tilde{\alpha}|^{\gamma+\sigma}.
\end{equation*}
Proceeding as above, one obtains the analogous bound for $Y_2$ and then from \eqref{H1H1}
$$|H_1-\mathcal{H}_1|\leq c(\|F(Z)\|_{L^\infty},\|\partial_{\alpha}Z\|_{C^{\gamma+\sigma}})|\alpha-\tilde{\alpha}|^{\gamma+\sigma}.$$
The same argument works for the difference $H_2-\mathcal{H}_2$, so that joining \eqref{G1H1H2} and \eqref{G1H1H22}, we can write
\begin{equation}\label{G}
\begin{aligned}
|G_1(\alpha,\beta)-G_1(\tilde{\alpha},\beta)-\mathcal{G}_1(\alpha,\beta)+\mathcal{G}_1(\tilde{\alpha},\beta)|&\leq |H_1-\mathcal{H}_1|+|H_2-\mathcal{H}_2|\\
&\leq c(\|F(Z)\|_{L^\infty},\|\partial_{\alpha}Z\|_{C^{\gamma+\sigma}})|\alpha-\tilde{\alpha}|^{\gamma+\sigma}.
\end{aligned}
\end{equation}
Since the term corresponding to $G_2$, $\mathcal{G}_2$ in \eqref{G1G2}, \eqref{G1G22} is completely analogous, the same bound \eqref{G} holds for $G$. Therefore, introducing this estimate into $P_8$ \eqref{O2}, we obtain that
\begin{equation}\label{P5}
\begin{aligned}
P_8&\leq c(\|F(Z)\|_{L^\infty},\|Z\|_{C^{1+\gamma+\sigma}},\eta)|\alpha-\tilde{\alpha}|^{\gamma+\sigma}(1-\log{|\alpha-\tilde{\alpha}|})\\
&\leq c(\|F(Z)\|_{L^\infty},\|Z\|_{C^{1+\gamma+\sigma}},\eta)|\alpha-\tilde{\alpha}|^{\gamma}.
\end{aligned}
\end{equation}
Joining the above bounds \eqref{P1P2},\eqref{P3}, \eqref{P4} and \eqref{P5} and going back to \eqref{O2} we find that
\begin{equation*}
|O_2(\alpha)-O_2(\tilde{\alpha})|\leq c(\|\theta_0\|_{L^\infty},\|F(Z)\|_{L^\infty},\|Z\|_{C^{1+\gamma+\sigma}})|h|^\gamma,
\end{equation*}
which concludes the subsection for $O_2$.
\vspace{0.2cm}

 This last bound combined with \eqref{O1bound} allow us to estimate \eqref{M1}
\begin{equation*}
|M_1|\leq c(\|\theta_0\|_{L^\infty},\|F(Z)\|_{L^\infty},\|Z\|_{C^{1+\gamma+\sigma}},\eta)|h|^\gamma,
\end{equation*}
which jointly to \eqref{M2bound} gives in \eqref{L6} that
\begin{equation*}
|L_6|\leq c(\|\theta_0\|_{L^\infty},\|F(Z)\|_{L^\infty},\|Z\|_{C^{1+\gamma+\sigma}},\eta)|h|^\gamma,
\end{equation*}
and thus the H\"older estimate of $I_1$ \eqref{I1I2I3} is concluded. Since we already have the bounds \eqref{I2I3bound}, formula \eqref{I1I2I3decomp} shows that the proof is ended.

\qed


\subsection*{{\bf Acknowledgments}}
This research was partially supported by the grant MTM2014-59488-P (Spain) and by the ERC through the Starting Grant project H2020-EU.1.1.-639227.
EGJ was supported by MECD FPU grant from the Spanish Government.

\vspace{-0.2cm}

\begin{quote}
\begin{tabular}{ll}
\textbf{Francisco Gancedo}\\
{\small Departamento de An\'{a}lisis Matem\'{a}tico $\&$ IMUS}\\
{\small Universidad de Sevilla}\\
{\small C/ Tarfia s/n, Campus Reina Mercedes, 41012 Sevilla, Spain}\\
{\small Email: fgancedo@us.es}
\end{tabular}
\end{quote}

\begin{quote}
\begin{tabular}{ll}
\textbf{Eduardo Garc\'ia-Ju\'arez}\\
\textbf{Former Address}\\
{\small Departamento de An\'{a}lisis Matem\'{a}tico $\&$ IMUS}\\
{\small Universidad de Sevilla}\\
{\small C/ Tarfia s/n, Campus Reina Mercedes, 41012 Sevilla, Spain}\\
\textbf{Current Address}\\
{\small Department of Mathematics}\\
{\small University of Pennsylvania}\\
{\small David Rittenhouse Lab., 209 South 33rd St., Philadelphia, PA 19104, USA}\\
{\small Email: edugar@math.upenn.edu}
\end{tabular}
\end{quote}

\end{document}